\documentclass{article}
\usepackage[english]{babel}
\usepackage[utf8]{inputenc}
\usepackage{amsmath}
\usepackage{graphicx}
\usepackage{amsfonts}
\usepackage{enumitem}
\usepackage{amssymb}
\usepackage[colorinlistoftodos]{todonotes}
\usepackage{geometry}
\usepackage{amsthm}
\usepackage{comment}
\usepackage{hyperref}
\usepackage[initials,nobysame]{amsrefs}
\usepackage{cleveref}
\usepackage{tikz}
\usetikzlibrary{math, calc}

\newcommand{\abs}[1]{\left\lvert#1\right\rvert}

\newcommand{\norm}[1]{\left\lVert#1\right\rVert}

\newcommand{\floor}[1]{\left\lfloor #1 \right\rfloor}
\newcommand{\ceil}[1]{\left\lceil #1 \right\rceil}
\newcommand{\paren}[1]{\left( #1 \right)}
\newcommand{\sqb}[1]{\left[ #1 \right]}

\newcommand{\setcond}[2]{\left\{ #1 \;\middle\vert\; #2 \right\}}

\newcommand{\one}{\mathbb{I}}

\newcommand{\pl}{\partial}

\newcommand{\on}{\operatorname}

\newcommand{\CC}{\mathbb{C}}

\newcommand{\EE}{\mathbb{E}}

\newcommand{\RR}{\mathbb{R}}
\newcommand{\PP}{\mathbb{P}}
\newcommand{\NN}{\mathbb{N}}
\newcommand{\ZZ}{\mathbb{Z}}
\newcommand{\QQ}{\mathbb{Q}}

\newcommand{\cB}{\mathcal{B}}

\newcommand{\cP}{\mathcal{P}}

\newcommand{\cS}{\mathcal{S}}

\newtheorem{thm}{Theorem}[section]
\newtheorem{lem}[thm]{Lemma}
\newtheorem{cor}[thm]{Corollary}

\newtheorem{prop}[thm]{Proposition}
\newtheorem{prob}[thm]{Problem}
\newtheorem*{rmk}{Remark}

\theoremstyle{definition}
\newtheorem{defn}[thm]{Definition}

\crefname{thm}{theorem}{theorems}
\Crefname{thm}{Theorem}{Theorems}
\crefname{lem}{lemma}{lemmas}
\Crefname{lem}{Lemma}{Lemmas}
\crefname{cor}{corollary}{corollaries}
\Crefname{cor}{Corollary}{Corollaries}
\crefname{conj}{conjecture}{conjectures}
\Crefname{conj}{Conjecture}{Conjectures}
\crefname{claim}{claim}{claims}
\Crefname{claim}{Claim}{Claims}
\crefname{obs}{observation}{observations}
\Crefname{obs}{Observation}{Observations}
\crefname{prop}{proposition}{propositions}
\Crefname{prop}{Proposition}{Propositions}
\crefname{prob}{problem}{problems}
\Crefname{prob}{Problem}{Problems}
\crefname{defn}{definition}{definitions}
\Crefname{defn}{Definition}{Definitions}
\crefname{example}{example}{examples}
\Crefname{example}{Example}{Examples}

\AddToHook{env/thm/begin}{\crefalias{thm}{thm}}
\AddToHook{env/lem/begin}{\crefalias{thm}{lem}}
\AddToHook{env/cor/begin}{\crefalias{thm}{cor}}
\AddToHook{env/conj/begin}{\crefalias{thm}{conj}}
\AddToHook{env/claim/begin}{\crefalias{thm}{claim}}
\AddToHook{env/obs/begin}{\crefalias{thm}{obs}}
\AddToHook{env/prop/begin}{\crefalias{thm}{prop}}
\AddToHook{env/prob/begin}{\crefalias{thm}{prob}}
\AddToHook{env/defn/begin}{\crefalias{thm}{defn}}
\AddToHook{env/example/begin}{\crefalias{thm}{example}}

\DeclareMathOperator{\lcm}{lcm}
\DeclareMathOperator{\Vol}{Vol}
\DeclareMathOperator{\Var}{Var}

\DeclareMathOperator{\rank}{rank}
\newcommand{\ord}{\mathrm{ord}}

\title{Efficient reconfiguration, polynomial decompositions and approximate counting}
\author{Peter Keevash\thanks{Mathematical Institute, University of Oxford, UK. Email: {\tt keevash@maths.ox.ac.uk}. Research supported by ERC Advanced Grant 883810.} \and Jeck Lim\thanks{Institute for Applied Mathematics, University of Bonn, 53115 Bonn, Germany. Email: \texttt{jlim@uni-bonn.de}. Supported by the Deutsche Forschungsgemeinschaft (DFG, German Research Foundation) under Germany's Excellence Strategy -- EXC-2047/2 -- 390685813. Part of this work was done while the author was at the University of Oxford supported by ERC Advanced Grant 883810.}}
\date{}

\begin{document}
\maketitle

\begin{abstract}
    We consider reconfigurations of integer assignments to planar lattice regions that preserve specified line totals,
    which are motivated by questions of combinatorial construction, algorithmic sampling and approximate counting.
    An elementary prototype is provided by contingency tables,
    i.e.~two-dimensional grids of non-negative integers with prescribed row and column sums,
    which can be reconfigured by basic switches of the form
    $\left(\begin{smallmatrix}1&-1\\-1&1\end{smallmatrix}\right)$.
    More difficult combinatorial instances include reconfiguration
    of Latin squares or queens configurations.
    Our main technical ingredient is the following polynomial decomposition lemma.
    Fix a positive integer $k$. Let $F(x)\in \mathbb{Z}[x]$ be a polynomial of degree at most $n$
    such that $(1-x)^k\mid F(x)$ and the sum of the absolute values of the coefficients of $F$ is $m$.
    We prove that $F$ can be written as a sum of $O(m+n)$ basic polynomials of degree at most $n$ of the form
    $\pm x^a(1-x^{a_1})\cdots (1-x^{a_k})$ where $a,a_1,\ldots,a_k$ are non-negative integers.
    The proof uses the existence of a miraculous high-dimensional coupling
    of signed distributions on cross polytopes that has independent interest.
    We will apply this lemma to questions of combinatorial construction in other papers;
    in this paper we illustrate its application to approximate counting.
\end{abstract}


\section{Introduction}

Consider the set $S$ of non-negative integer matrices $A = (A_{ij})_{i,j=1}^n$
with prescribed row sums $r_i = \sum_j A_{ij}$ and column sums $c_j = \sum_i A_{ij}$.
Can we approximate $|S|$? Efficiently sample from $S$? Efficiently navigate between elements of $S$?
These questions are classical in combinatorics \cite{Ryser1957,Barvinok2009,Barvinok2010,BarvinokHartigan2010},
random graph theory \cite{Greenhill2021,CanfieldGreenhillMcKay2008},
computer science \cite{JerrumValiantVazirani1986}
and statistics \cite{DiaconisSturmfels1998,AlmendraDeLoeraPetrovic2024}.
Switchings provide a simple and powerful proof technique, based on the observations
that any difference $D=A-B$ with $A,B\in S$ has all row and column sums zero,
and so can be efficiently decomposed into summands, each of which
is a basic switch of the form $\left(\begin{smallmatrix}1&-1\\-1&1\end{smallmatrix}\right)$
for some choice of two rows and two columns. It is not hard to see that for any such nonzero $D$
with all row and column sums zero there is a basic switch that decreases $|D|:=\sum_{i,j} |D_{ij}|$
by at least $2$, so the required number of summands is $O(|D|)$.

Going up one dimension, a Latin square of order $n$, or $3$-dimensional permutation,
can be represented by a $\{0,1\}$-valued  $3$-tensor $A = (A_{ijk})_{i,j,k=1}^n$,
where the indices represent rows, columns and symbols, with the property that
every axis-parallel line contains exactly one $1$.
The corresponding basic switches are $2 \times 2 \times 2$ tensors
obtained by stacking a two-dimensional basic switch and its negation;
see \cite{JacobsonMatthews1996,Cavenagh2008} for the theory
of navigation via switches for Latin squares.

Similar complexity can be seen in two dimensions when more directions are constrained,
such as in the approximate counting problem for queens configurations \cite{BowtellKeevash2021,Simkin2023},
via the designs counting strategy developed in \cite{keevash2018counting}, which has two stages:
(i) a randomised construction provides the correct entropy, but only gives approximate objects, then
(ii) an absorption procedure with small entropy cost transforms approximate objects into exact objects.
Our paper is motivated by the question of implementing (ii) for more general systems of linear constraints.
Thus it is similar in spirit to \cite{KeevashLong2020}, but in an opposite regime: instead of high-dimensional cubes
with few constraints, we will consider two-dimensional regions with many constraints, obtaining absorption and
approximate counting results that do not follow from the general machinery in \cite{Barvinok2009,Barvinok2010,BarvinokHartigan2010}.

To state our approximate counting result we need the following notation.
Let $v_1,\dots,v_k \in \ZZ^2$ be non-parallel primitive vectors, fix $c_1,\dots,c_k \in \NN$
and consider a zonotope $Z = \sum_{i=1}^k [-c_i,c_i] v_i$
that is the Minkowski sum of segments in these directions.
Let $\Omega_n = nZ \cap \ZZ^2$ and let $A_n$ be the incidence matrix
of all lines parallel to some $v_i$ that meet $\Omega_n$,
i.e.~$(A_n)_{\ell,x} = 1_{x \in \ell}$ for each such $\ell$ and $x \in \Omega_n$.
Given $b \in A_n \{0,1\}^{\Omega_n}$ let $N(b)$ count $x \in \{0,1\}^{\Omega_n}$ with $A_n x = b$.
Choose $p \in [0,1]^{\Omega_n}$ with $A_n p = b$
to maximise its entropy $H(p)$.

\begin{thm} \label{thm:count}
    In the above setting, for any $c>0$ and $n$ sufficiently large,
    if $p_v \in [c,1-c]$ for all $v \in \Omega_n$ then $N(b) = 2^{H(p)+o(n^{1.51})}$.
\end{thm}

\subsection{Absorption via polynomial decomposition} \label{subsec:absorb}

We will adopt an algebraic perspective formulated in \cite{HajduTijdeman2001} and further developed in \cite{StolkBatenburg2010}
in the context of discrete tomography \cite{gardner1997discrete,gardner2006geometric}, a field of research motivated
by a practical question posed in 1994 by Larry Shepp of AT\&T Bell Labs: how can one reconstruct a finite subset $F$ of a lattice,
given a collection of `X-rays' that count points of $F$ along some given directions?
(This problem is NP-complete \cite{gardner1999computational},
so typically one seeks approximate solutions or makes further assumptions, e.g. convexity.)
We represent $f:\ZZ^2 \to \ZZ$ of finite support by its two-variable generating function
$p_f(X,Y) = \sum_{s,t} f(s,t) X^s Y^t$, which is a Laurent polynomial in $R=\ZZ[X,X^{-1},Y,Y^{-1}]$.
X-rays are extracted by appropriate substitutions, e.g. $p_f(z,1) = \sum_{s,t} f(s,t) z^s$
is the one-variable generating function for line sums in the $Y$-direction.
In general, $p_f(z^b,z^{-a})$ is the generating function for line sums
in some direction $(a,b) \in \ZZ^2$ with $a,b$ coprime:
the coefficient of $z^c$ sums all $f(s,t)$ with $bs-at=c$.
This leads to the characterisation of `switching components' established in \cite{HajduTijdeman2001}:
given some constraint directions $v_i = (a_i,b_i) \in  \ZZ^2$ for $i=1,\dots,k$,
a function $f:\ZZ^2 \to \ZZ$ has zero sums in all directions parallel to any $v_i$
precisely when all projections $\pi_{v_i}(f) := p_f(z^{b_i},z^{-a_i})=0$;
equivalently, $p_f$ is divisible in $R$ by $\prod_{i=1}^k (X^{a_i} - Y^{-b_i})$.

To motivate our absorption lemma below, consider the problem of placing generalised `queens'
so that each line $\ell$ parallel to any $v_i$ contains precisely $n_\ell$ queens,
for some fixed choice of the $(n_\ell)$ where finitely many are nonzero.
We relax our task by allowing queens with (positive or negative) integer multiplicities
(which can be implemented via suitable gadgets in the absorption strategy mentioned above).
We consider each direction $v_i$ in turn and resolve all lines parallel to $v_i$ without breaking any previously resolved lines.
For $v_1$ this is trivial: we place $n_\ell$ queens arbitrarily on each $\ell$ parallel to $v_1$.
We then update each $n_\ell$ to its remaining requirement, so now $n_\ell=0$ for each $\ell$ parallel to $v_1$.
For $v_2$, we need to find $f:\ZZ^2 \to \ZZ$ with $\pi_{v_1}(f)=0$ and $\pi_{v_2}(f)=h_2(z)$,
where $h_2(z)$ is the generating function for required line sums in the $v_2$ direction.
This requires $X^{a_1}-Y^{-b_1} \mid p_f$, so $h_2(z)$ must be divisible in $\ZZ[z,z^{-1}]$
by $\pi_{v_2}(X^{a_1}-Y^{-b_1}) = z^{b_2 a_1} - z^{a_2 b_1}$,
or equivalently by $z^{m_{12}}-1$ where $m_{ij}=|\det(v_i,v_j)|$.
Given such $h_2(z)$, one can greedily construct $f$ by placing `cancelling' pairs of queens,
each having one positive and one negative queen on some line parallel to $v_1$.
Under $\pi_{v_2}$ such a pair maps to a polynomial of the form $z^c (z^{dm_{12}}-1)$ with $c,d \in \ZZ$.
Thus we can repeatedly eliminate nonzero coefficients in $h_2(z)$
with exponents differing by a multiple of $m_{12}$ until $h_2(z)$ becomes the zero polynomial.
We then update each $n_\ell$ to its remaining requirement,
so now $n_\ell=0$ for each $\ell$ parallel to $v_1$ or $v_2$.

For $v_3$, we need to find $f:\ZZ^2 \to \ZZ$ with $\pi_{v_1}(f)=\pi_{v_2}(f)=0$ and $\pi_{v_3}(f)=h_3(z)$,
where $h_3(z)$ is the generating function for required line sums in the $v_3$ direction.
As above, $h_3(z)$ must be divisible by $(z^{m_{13}}-1)(z^{m_{23}}-1)$.
To construct such $h_3(z)$, our basic switches must go up one dimension: instead of cancelling pairs,
we use the $\left(\begin{smallmatrix}1&-1\\-1&1\end{smallmatrix}\right)$ switches discussed above,
which here means placing positive queens at $v,v+cv_1+dv_2$ and negative queens at $v+cv_1,v+dv_2$
for some $v \in \ZZ^2$ and $c,d \in \ZZ$. Under $\pi_{v_3}$ such a switch maps to a polynomial
of the form $z^c (z^{d_1 m_{13}}-1) (z^{d_2 m_{23}}-1)$ with $c,d_1,d_2 \in \ZZ$.
The switching procedure mentioned above for contingency tables
can be adapted to provide the required decomposition of $h_3(z)$ into switches.

For $v_{i+1}$ with $i \ge 3$, the required choice of $i$-dimensional switches translates into the problem of decomposing
some $h_{i+1}(z)$ divisible by $\prod_{j\le i} (z^{m_{j,i+1}}-1)$ into basic polynomials of the form $z^c \prod_{j\le i}(  z^{d_j m_{j,i+1}}-1)$.
This is a substantially more difficult problem, as will be explained below.
It is already apparent that a naive greedy argument fails,
as if $F$ is supported on a random set of size $k+1$
then no basic switch can have support close to that of $F$.
Moreover, our lower bounds on filling norms mentioned below
imply that no greedy algorithm
analogous to those in small dimensions can work in higher dimensions,
so more sophisticated ideas are required in our solution.

\subsection{Polynomial decompositions and filling norms}

To state our main result, we require some definitions.
Let $m = (m_1,\dots,m_k) \in \ZZ_+^k$.
We call a polynomial \emph{$m$-basic} if it has the form
$\pm x^a(1-x^{m_1a_1})\cdots (1-x^{m_ka_k})$
with $a,a_1,\ldots,a_k \in \NN$.
If all $m_i=1$ we simply say \emph{basic}, with $k$ understood from the context.
For a polynomial $P=\sum_i a_i x^i \in \ZZ[x]$, or more generally,
a Laurent polynomial $P\in \ZZ[x,x^{-1}]$,
we define $|P|:=\sum_i |a_i|$ to be the size of $P$.

The following theorem provides the quantitative decomposition result
required to implement the absorption strategy outlined in the previous subsection.

\begin{thm} \label{thm:decomp-rect}
    Fix $m=(m_1,\ldots,m_k) \in \ZZ_+^k$.
    There exists $C>0$ such that the following holds.
    Let $F\in \ZZ[x]$ be a polynomial of degree at most $n$ such that $(1-x^{m_1})\cdots (1-x^{m_k})\mid F$.
    Then $F$ can be written as a sum of at most $C\max(|F|,n)$ polynomials,
    each of which is $m$-basic and of degree at most $n$.
\end{thm}

It is instructive to interpret this result in Gromov's language of filling norms \cite[Section 2.3]{gromov2010singularities},
which is foundational for the theory of high-dimensional expanders \cite{lubotzky2018high}.
Given a homomorphism $\pl: A \to B$ between two normed abelian groups,
we define the \emph{filling norm} of $b \in B$ by $\|b\|_{\on{fill}} = \inf_{a \in \pl^{-1}(b)} \|a\|$.
The  \emph{filling profile} at $\beta$ is $\sup_{\|b\|=\beta}  \|b\|_{\on{fill}}/\|b\|$.
In the setting of \Cref{thm:decomp-rect}, we consider $B = \ZZ[x]$ and $A = \ZZ^{\cB}$,
where $\cB$ is the set of $m$-basic polynomials and $\pl(a) = \sum_{p \in \cB} a_p p \in B$.
Here we have temporarily suppressed the degree conditions: one might think these
are a technical distraction (they are needed for combinatorial applications)
and that the filling profile is uniformly bounded.
However, we will show in Section \ref{subsec:lower}
that this is not the case: for $k=3$ one can find $F$ of degree at most $n$
with $|F|=6$ such that $\|F\|_{\on{fill}} = \Omega(\log n)$. This phenomenon disappears when
we move to real coefficients, when the filling profile is uniformly bounded,
and furthermore we can construct a contraction (or cone) in the sense of  \cite[Section 2.6]{gromov2010singularities}.

\begin{thm} \label{thm:real}
    Fix $m=(m_1,\ldots,m_k) \in \ZZ_+^k$ and $n \in \ZZ_+$.
    Let $B$ be the set of polynomials in $\RR[x]$ of degree at most $n$.
    Let $\cB$ be the set of $m$-basic polynomials of degree at most $n$.
    Let $A = \RR^{\cB}$ and define $\pl: A \to B$ by $\pl(a) = \sum_{p \in \cB} a_p p$.
    Then there is a homomorphism $\delta:B \to A$ such that  for all $b \in \pl(A)$ we have
    $\pl \circ \delta(b)=b$ and $|\delta(b)| \le C|b|$ for some $C$ depending only on $m$.
\end{thm}

This sensitivity to the choice of coefficients is both geometrically interesting and significant for our proof,
as we use (an approximate variant of) \Cref{thm:real} in the proof of \Cref{thm:decomp-rect}.
This follows an iterative absorption scheme introduced in version 1 of \cite{keevash2014existence},
which reduces $F\in \ZZ[x]$ as in  \Cref{thm:decomp-rect} to zero by alternating two steps:
(i) cancellation steps that significantly reduce $|F|$, and
(ii) focusing steps that reduce the degree of $F$ by a constant factor,
while not unduly increasing $|F|$. Here (ii) is straightforward and (i) needs the key idea:
we can convert real solutions to approximate integer solutions by randomised rounding.
The probabilistic viewpoint is suggestive for the proof of \Cref{thm:real}:
we will find approximate real solutions via a miraculous coupling discussed in the next subsection.

It will be helpful to first prove our polynomial decomposition results
in the setting of basic polynomials, i.e.~when all $m_i=1$.
This reduces the full theorem to an instance of the following problem,
which can be thought of as an instance of Gromov filling for the group $G=I/J$.

\begin{prob} \label{prob:general-ideal}
    Let $J\leq I$ be ideals of a polynomial ring and $\cS$ be a fixed class of `small' polynomials in $I$.
    Given $f\in I$, can we find a `small' sum $g$ of elements of $\cS$ such that $f-g\in J$?
\end{prob}

We will not study this problem systematically, but instead we will prove
two results of this type that we need for particular applications.
The first, \Cref{lem:add-cyclotomic}, is used to deduce \Cref{thm:decomp-rect} from the basic case.
The second, \Cref{prop:two-var}, is a technical two-variable lemma that we will use in \cite{KL26-2}
to provide a positive answer to the following question of Sagdeev \cite{MO}.

\begin{prob}
    Let $D_1,\ldots,D_k\subset \ZZ^2$ be `dominoes', i.e. two-element subsets of $\ZZ^2$.
    Is there a finite non-empty $B\subset \ZZ^2$ that for each $i$ can be partitioned into disjoint translates of $D_i$?
\end{prob}

\subsection{Motivation for the cancellation step}

Let $F=\sum_{i=1}^m \epsilon_i x^{d_i}\in \ZZ[x]$ be a polynomial of degree at most $n$ such that $(1-x)^k\mid F$, with $\epsilon_i\in \{\pm 1\}$ and $m=|F|$.
The basic idea of the cancelling step is to replace each monomial $x^d$ with $x^d-P_d$,
for some random independent basic polynomial $P_d$ containing the monomial $x^d$,
chosen so that the remainder is small in expectation,  and so with high probability by concentration.
We do not know of an easier way to make the remainder small than to make it zero, i.e.~
\begin{equation} \label{eq:cancel-expectation}
    \sum_i \EE[\epsilon_i (x^{d_i}-P_{d_i})]=0.
\end{equation}
The most naive sufficient condition for cancellation is to have $\EE[x^d-P_d]=0$.
However, this is impossible since the sum of the coefficients of $P_d$ is always zero, so the sum of the coefficients of $\EE[x^d-P_d]$ is always 1.
The next attempt would be to ensure
\begin{equation} \label{eq:const-coeff}
    \EE[x^d-P_d]=\frac{1}{n+1}(1+x+\cdots+x^n).
\end{equation}
Since $1-x\mid F$, we have $\sum_i \epsilon_i=0$, so summing gives (\ref{eq:cancel-expectation}).
This works well for $k=1$, but (\ref{eq:const-coeff}) is impossible for larger $k$ due to higher moment constraints.
Indeed, suppose $\EE[x^d-P_d]=\sum_{i=0}^n p(i) x^i$.
Since $(1-x)^k\mid P_d$, differentiating gives $P_d^{(j)}(1)=0$ for $j=0,\ldots,k-1$, i.e.
\begin{equation} \label{eq:moment-constraints}
    \sum_{i=0}^n p(i) i^j=d^j \qquad \text{for } j=0,\ldots,k-1.
\end{equation}
These moment constraints are not satisfied by (\ref{eq:const-coeff}) in general for $k\geq 2$.
Instead we can aim for
\begin{equation} \label{eq:poly-coeff}
    \EE[x^d-P_d]=\sum_{i=0}^n p_d(i) x^i,
\end{equation}
where $p_d$ is a polynomial of degree at most $k-1$. This will be enough to imply (\ref{eq:cancel-expectation}).
Indeed, summing (\ref{eq:poly-coeff}) over all the $d_i$'s gives
\(    \sum_{i=1}^m \EE[\epsilon_i (x^{d_i}-P_{d_i})]  =  \sum_{i=0}^n \tilde{p}(i) x^i \)
for some polynomial $\tilde{p}(t) = \sum_{j=0}^{k-1} c_j t^j$ of degree at most $k-1$.
Since $(1-x)^k\mid F$ we have $\sum_{i=0}^n i^j \tilde{p}(i)=0$ for $j=0,\ldots,k-1$.
Summing the $j$th equation multiplied by $c_j$ gives $\sum_{i=0}^n \tilde{p}(i)^2=0$,
and thus $\tilde{p}\equiv 0$, so (\ref{eq:cancel-expectation}) holds.
We will see that such a choice for $p_d$ satisfying (\ref{eq:moment-constraints}) is possible; in fact, the choice is unique.

The aim now is to construct a distribution on basic polynomials $P_d$ such that (\ref{eq:poly-coeff}) holds;
we will ensure the stronger condition that  each non-trivial monomial of $P_d$ (besides $x^d$)
has a distribution with a degree $k-1$ polynomial density function.
Also, we require $P_d = x^d(1-x^{a_1})\cdots (1-x^{a_k})$ of degree at most $n$,
so we can parameterise the support of $P_d$ by the polytope  $\cP_d$
consisting of all $(a_1,\ldots,a_k)\in \ZZ^k$ satisfying
\[
    0\leq d+a_S\leq n \quad \text{for all }S\subseteq [k],\qquad \text{where } a_S:=\sum_{i\in S} a_i.
\]
Thus we seek a distribution on $\cP_d$ such that the density function of each exponent $d+a_S$ is a polynomial of degree at most $k-1$, supported on $[0,n]$.

\begin{figure}[htbp]
    \centering
    \begin{tikzpicture}[x=4.8cm,y=4.8cm,font=\small,line join=round]
        \definecolor{pdDensity}{RGB}{41,97,143}
        \foreach \pdSliceOffset/\pdShift/\pdPanel in
            {0/0/{(a) $\cP_{n/2}$},0.5/1.25/{(b) $\cP_{n/4}$}} {
                \begin{scope}[shift={(\pdShift,0)}]
                    \foreach \pdStep in {160,159,...,0} {
                            \pgfmathsetmacro{\pdLevel}{\pdSliceOffset+(1-\pdSliceOffset)*\pdStep/160}
                            \pgfmathsetmacro{\pdTint}{100*\pdLevel}
                            \fill[pdDensity!\pdTint!white]
                            ({-\pdLevel/2},{-\pdSliceOffset*sqrt(3)/6}) --
                            ({(\pdSliceOffset-\pdLevel)/4},{(\pdSliceOffset+3*\pdLevel)*sqrt(3)/12}) --
                            ({(\pdLevel-\pdSliceOffset)/4},{(\pdSliceOffset+3*\pdLevel)*sqrt(3)/12}) --
                            ({\pdLevel/2},{-\pdSliceOffset*sqrt(3)/6}) --
                            ({(\pdSliceOffset+\pdLevel)/4},{(\pdSliceOffset-3*\pdLevel)*sqrt(3)/12}) --
                            ({-(\pdSliceOffset+\pdLevel)/4},{(\pdSliceOffset-3*\pdLevel)*sqrt(3)/12}) -- cycle;
                        }
                    \coordinate (pdA) at (-.5,{-\pdSliceOffset*sqrt(3)/6});
                    \coordinate (pdB) at ({(\pdSliceOffset-1)/4},{(\pdSliceOffset+3)*sqrt(3)/12});
                    \coordinate (pdC) at ({(1-\pdSliceOffset)/4},{(\pdSliceOffset+3)*sqrt(3)/12});
                    \coordinate (pdD) at (.5,{-\pdSliceOffset*sqrt(3)/6});
                    \coordinate (pdE) at ({(\pdSliceOffset+1)/4},{(\pdSliceOffset-3)*sqrt(3)/12});
                    \coordinate (pdF) at ({-(\pdSliceOffset+1)/4},{(\pdSliceOffset-3)*sqrt(3)/12});
                    \draw[black!55,line width=.4pt]
                    (pdA)--(pdD) (pdB)--(pdE) (pdC)--(pdF);
                    \ifdim\pdSliceOffset pt>0pt
                        \draw[pdDensity!65!black,line width=.7pt]
                        ({-\pdSliceOffset/2},{-\pdSliceOffset*sqrt(3)/6}) --
                        (0,{\pdSliceOffset*sqrt(3)/3}) --
                        ({\pdSliceOffset/2},{-\pdSliceOffset*sqrt(3)/6}) -- cycle;
                        \node[inner sep=1pt] at (0,0) {core};
                    \else
                        \fill[black!70] (0,0) circle[radius=.65pt];
                    \fi
                    \draw[pdDensity!65!black,line width=.8pt]
                    (pdA)--(pdB)--(pdC)--(pdD)--(pdE)--(pdF)--cycle;
                    \node[anchor=north,inner sep=2pt] at (0,-.54) {\pdPanel};
                \end{scope}
            }
        \shade[bottom color=white,top color=pdDensity] (1.98,-.38) rectangle (2.04,.38);
        \draw[black!60,line width=.35pt] (1.98,-.38) rectangle (2.04,.38);
        \foreach \pdValue/\pdLabel in {0/0,.5/{\tfrac12},1/1} {
                \draw[black!60,line width=.35pt]
                (2.04,{-.38+.76*\pdValue}) -- (2.06,{-.38+.76*\pdValue});
                \node[anchor=west,inner sep=2pt] at (2.06,{-.38+.76*\pdValue}) {$\pdLabel$};
            }
        \node[anchor=south,inner sep=2pt] at (2.01,.40) {$\rho/n$};
    \end{tikzpicture}
    \caption{Affine images of $\cP_d$ for $k=2$: (a) the midpoint $d=n/2$; (b) the slice $d=n/4$. The colour shows $\rho/n$, where $\rho=|n-2d-a_1-a_2|+|a_1|+|a_2|$; the density is proportional to this quantity.}
    \label{fig:hexagon-densities}
\end{figure}
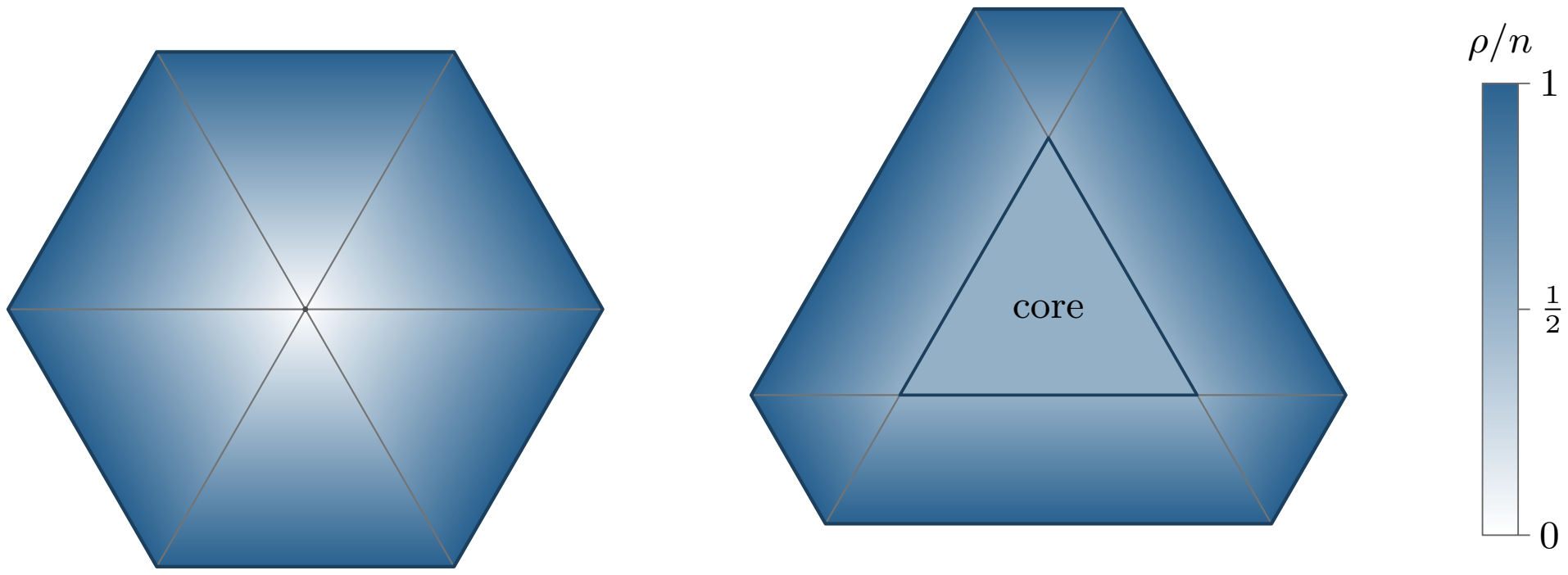

This might seem too much to ask for, but small cases show some promise.
Consider the case $k=2$. For the endpoints $d=0,n$, $\cP_d$ is a simplex, so the uniform distribution works. Another special case is the midpoint $d=n/2$. In this case, $\cP_{n/2}$ is an affine transformation of the regular hexagon. It is convenient to take $n\to\infty$ and consider the continuous limit instead. Then we can find a distribution on the regular hexagon such that the projections along the three directions of the hexagon are uniform on the support. Indeed, we can put the Archimedean measure on the inscribed disk (this is the projection of the uniform surface measure on the sphere, see for example \cite{angel2007random}). Its density becomes unbounded at the boundary of the disk, which is detrimental for our purposes. Instead, experimenting with small $n$ reveals a distribution with bounded, piecewise linear density and uniform marginals along the three directions. Its density starts from 0 at the center and grows linearly to the boundary of the hexagon. This distribution can be generalised to other $d$. The hexagon $\cP_d$ can be partitioned by its diagonals, forming a triangular core. The distribution is constant on the triangular core, and grows linearly to the boundary of $\cP_d$. Its three marginals are linear on the support. See \Cref{fig:hexagon-densities}.

Before we can generalise this to larger $k$, we make the following crucial observation.
The polytopes $\cP_d$ can be stacked together to form the cross polytope $\cP$ in $\RR^{k+1}$, defined by
\[
    \cP=\setcond{(a_0,a_1,\ldots,a_k)\in \RR^{k+1}}{\sum_i |a_i|\leq n}.
\]
Introducing the extra coordinate
\(
a_0=n-2d-\sum_{i=1}^k a_i
\)
identifies $\cP_d$ with the slice of $\cP$ at $\sum_i a_i=n-2d$. For $k=2$ the above distributions also stack nicely.
The resulting density function $Q$ on $\cP$ is linear in terms of the $\ell^1$ radius $\rho(x)=\sum_i |x_i|$, with $Q(0)=0$.

For larger $k$, we guess that such a density function $Q$ on $\cP$ exists, and is a function of the radius $\rho(x)$. Miraculously, this turns out to be true, and $Q$ is a polynomial of degree $k-1$ in $\rho$. After computing this polynomial for small $k$, we discovered that $Q$ is (proportional to) the Legendre polynomial. One caveat is that the Legendre polynomial takes negative values, so $Q$ is technically not a density function. Nonetheless, signed distributions are sufficient for our purposes. Given this explicit distribution, verifying the marginals is not difficult, despite being lengthy.

\subsection{A miraculous distribution}

For each integer $n\geq 0$, consider the Rodrigues formula
\[L_n(x):=\frac{1}{2^n n!}\frac{d^n}{dx^n}(x^2-1)^n\]
for the Legendre polynomial of degree $n$. We have the following properties of $L_n$.
\begin{enumerate}
    \item (Normalisation) $L_n(1)=1$ for all $n\geq 0$.
    \item (Orthogonality) For all $m,n\geq 0$, we have
    \[
        \int_{-1}^1 L_m(x)L_n(x) dx=\begin{cases}
            0,              & m\neq n \\
            \frac{2}{2n+1}, & m=n.
        \end{cases}
    \]
    \item (Parity) For all $n\geq 0$, we have $L_n(-x)=(-1)^n L_n(x)$.
    \item (Murphy formula) For all $n\geq 0$, we have
    \[L_n(x)=\frac{1}{2^n}\sum_{i=0}^n \binom{n}{i}^2 (x-1)^i (x+1)^{n-i}.\]
\end{enumerate}

We define the following explicit distribution on the cross-polytope $\setcond{x\in \RR^{k+1}}{\sum_{i=0}^k |x_i|\leq 1}$.

\begin{thm} \label{thm:crosspoly}
    Let $k$ be a positive integer and $Q:\mathbb{R}^{k+1}\to \mathbb{R}$ be given by $Q(x_0,\ldots,x_k)=L_{k-1}(|x_0|+\cdots+|x_k|)$ when $0\le |x_0|+\cdots+|x_k|\le 1$, and $Q(x)=0$ elsewhere. Let $l<k$ be a non-negative integer. For $r,s\in \mathbb{R}$, define $R_{l}(r,s)=\int_{V_{r,s}} Q(x_0,\ldots,x_k)\,dx$, where $V_{r,s}$ is the affine subspace given by $x_0+\cdots+x_k=r$ and $x_0+\cdots+x_l-x_{l+1}-\cdots-x_k=s$. Then the following hold.
    \begin{enumerate}
        \item $R_l(r,s)$ is a polynomial in the domain $(-1,1)^2$.
        \item $\deg_r R_l=\deg_s R_l\leq k-1$.
        \item $\int_{-1}^1 R_l(r,s)\, ds=2/k!$ for each $-1<r<1$.
    \end{enumerate}
\end{thm}

\begin{rmk}
    In the statement and the proof, a hyperplane of the form $\pm x_1\pm\cdots\pm x_l=c$ is equipped with the measure obtained by eliminating any one coordinate and using standard Lebesgue measure on the remaining $l-1$ coordinates. On $V_{r,s}$, we use the product of these measures on the two blocks $x_0+\cdots+x_l=(r+s)/2$ and $x_{l+1}+\cdots+x_k=(r-s)/2$.
\end{rmk}

Let us pause to appreciate the miracles required for Theorem~\ref{thm:crosspoly} to hold. The function $Q$ is piecewise polynomial of degree $k-1$. Suppose $L_{k-1}$ were replaced by a generic polynomial of degree $k-1$. Since $R_l$ is a $(k-1)$-dimensional integral of $Q$, we expect $R_l$ to be piecewise polynomial of degree $2k-2$. There are four pieces, cut out by the lines $r+s=0$ and $r-s=0$. By symmetry, it is not hard to see that $R_l(r,s)=R_l(s,r)=R_l(-s,-r)$. For $R_l$ to be a global polynomial, its polynomial must satisfy $R_l(r,s)=R_l(s,r)=R_l(-s,-r)$. This means that the coefficient of some monomial $r^a s^b$ must equal the coefficient of $r^b s^a$, and is zero if $a+b$ is odd. This is a total of $3k(k-1)/2$ linear constraints on the coefficients of $R_l$, which is unreasonable to expect since we only have $k$ degrees of freedom in choosing the coefficients of $L_{k-1}$. Furthermore, all this is only for satisfying (1) for a single $l$! It is true that if (1) holds for some $l$, then it also holds for $k-l-1$. However, there is no reason to expect (1) to also hold for other values of $l$.
Later, we will see that $\deg_s R_l\leq k-1$ on one of the pieces for a generic polynomial. Thus (1) implies symmetry, giving (2). As far as we can tell, property (3) does not immediately follow from (1) and (2).

\subsection{Organisation of the paper}

In the next section we prove \Cref{thm:crosspoly} regarding the Legendre distribution.
We apply it in Section 3 to prove \Cref{thm:real} and then \Cref{thm:decomp-rect}
in the case of basic polynomials. Then in Section 4 we prove these results in full generality
via a solution to \Cref{prob:general-ideal} for certain ideals defined by cyclotomic polynomials.
Section 5 contains the proof of \Cref{prop:two-var}.
We conclude in Section \ref{sec:count} with the proof of the approximate counting result \Cref{thm:count}.

\subsection*{Statement on AI use}
A complete proof of the main theorems (\Cref{thm:decomp-rect} and all its prerequisites) was written by the authors without significant AI assistance by October 2025.
After the proofs were written, GPT-6 Astra/6.1 Sol was used for proofreading, editing and figure generation;
it also provided the simple proof of the filling norm lower bound presented in Section \ref{subsec:lower}.

\section{Analysis of the Legendre distribution}

Here we verify the miraculous properties of the Legendre distribution recorded in \Cref{thm:crosspoly}.
Our proof is a direct computation; we do not know a conceptual reason for these miracles.

\subsection{Polytope slices}

We first compute certain quantities related to slices of the $n$-dimensional cross-polytope
\[
    \cP_n=\setcond{x\in \RR^n}{\sum_{i=1}^n |x_i|\leq 1},
\]
which has volume $\Vol_n(\cP_n)=2^n/n!$.
For a non-negative integer $n$, let
\[
    f_n(a,b)=\frac{1}{2^n}\sum_{i=0}^n \binom{n}{i}^2 a^i b^{n-i},
\]
so that $L_n(x)=f_n(x-1,x+1)$.

Let $V_s$ be the affine hyperplane $\sum_{i=1}^n x_i=s$. We shall first compute the slice volume $\Vol_{n-1}(\cP_n \cap V_s)$ for $-1\leq s\leq 1$.

Project $\cP_n \cap V_s$ onto $\RR^{n-1}$ by forgetting the last coordinate. We split this projected region into $2^{n-1}$ parts, corresponding to the sign pattern of $(x_1,\ldots,x_{n-1})$. The part with $l\geq 0$ positive signs, say $x_1,\ldots,x_l\geq 0$ and $x_{l+1},\ldots,x_{n-1}\leq 0$, consists of the points $(x_1,\ldots,x_{n-1})$ with $x_1+\ldots+x_l-x_{l+1}-\cdots-x_{n-1}+|s-x_1-\cdots-x_{n-1}|\leq 1$, or equivalently,
\begin{align*}
    x_1+\cdots+x_l            & \leq \frac{1+s}{2} \quad\text{and} \\
    -(x_{l+1}+\cdots+x_{n-1}) & \leq \frac{1-s}{2}.
\end{align*}
The volume of this set is exactly
\begin{align*}
    \paren{\frac{1+s}{2}}^l\frac{1}{l!}\cdot \paren{\frac{1-s}{2}}^{n-1-l}\frac{1}{(n-l-1)!} & = \frac{1}{2^{n-1}l!(n-1-l)!}(1+s)^l(1-s)^{n-l-1}.
\end{align*}
Summing over all the sign patterns, we obtain
\begin{align*}
    \Vol_{n-1}(\cP_n \cap V_s) & = \sum_{l=0}^{n-1} \binom{n-1}{l}\frac{1}{2^{n-1}l!(n-1-l)!}(1+s)^l(1-s)^{n-l-1} \\
                               & =\frac{1}{(n-1)!}f_{n-1}(1+s,1-s) \\
                               & =\frac{s^{n-1}}{(n-1)!}f_{n-1}(s^{-1}+1,s^{-1}-1) \\
                               & = \frac{s^{n-1}}{(n-1)!}L_{n-1}(s^{-1}).
\end{align*}
Expressions involving $s^{-1}$ are understood at $s=0$ by taking limits. For instance, the above expression is a polynomial in $s$ after expansion, so there is no problem taking $s=0$.

Next, define the $\ell^1$ radius $\rho(x):=\sum_{i=1}^n |x_i|$, and let $\nu_{n,s}$ be the pushforward under $\rho$ of the measure on $V_s$ specified above. So for $a\geq |s|$, we have
\[
    \nu_{n,s}([0,a])=\Vol_{n-1}(a\cP_n\cap V_s)
    =\frac{s^{n-1}}{(n-1)!}L_{n-1}(a/s),
\]
and $\nu_{n,s}([0,a])=0$ for $0\leq a<|s|$. Thus, for $n\geq 2$ and $a>|s|$, the density of $\nu_{n,s}$ is
\[
    S_{n-1}(s,a):=\frac{d}{da}\Vol_{n-1}(a\cP_n\cap V_s)
    =\frac{s^{n-2}}{(n-1)!}L_{n-1}'(a/s).
\]
At the critical radius $a=|s|$, the slice $a\cP_n\cap V_s$ is an $(n-1)$-simplex, and $\nu_{n,s}$ has a point mass at this radius of size
\[
    S_{n-1}(|s|):=\frac{|s|^{n-1}}{(n-1)!}.
\]
For $n=1$, the slice is a single point, so we set $S_0(s,a)=0$ and $S_0(a)=1$. In all cases, writing $\delta_t$ for a unit point mass at $t$, we have
\[
    \nu_{n,s}(da)=S_{n-1}(|s|)\delta_{|s|}(da)
    +\one(a>|s|)S_{n-1}(s,a)\,da.
\]
These formulas define polynomials $S_n(s,a)$ and $S_n(a)$, which we extend to all real arguments. We summarise their definitions and some properties which are easy to verify.
\begin{enumerate}
    \item $S_n(s,a)=\frac{s^{n-1}}{n!}L_n'(a/s)$ for $n\geq 1$, and $S_0(s,a)=0$.
    \item $S_n(a)=\frac{a^n}{n!}$.
    \item For $n\geq 1$, $S_{n}(s,a)$ is a homogeneous polynomial of degree $n-1$.
    \item $S_{n}(-s,a)=S_{n}(s,a)=(-1)^{n-1}S_n(s,-a)$.
    \item $S_{n}(a)$ is a polynomial of degree $n$ with $S_n(-a)=(-1)^n S_n(a)$.
\end{enumerate}

\subsection{Direct computations}

In this subsection, we prove \Cref{thm:crosspoly} by a (long) direct computation. Fix $l$ and write $R=R_l$ for simplicity. Fix $r,s$ and set $a=\frac{r+s}{2},b=\frac{r-s}{2}$. Then the subspace $V_{r,s}$ is defined by $x_0+\ldots+x_l=a$ and $x_{l+1}+\cdots+x_k=b$. We will assume for now that $1\geq r\geq |s|$, so that $a,b\geq 0$. Then we have
\begin{align*}
    R(r,s) & = \int_{V_{r,s}}Q(x)\ dx=\int_{V_{r,s}}L_{k-1}(|x_0|+\cdots+|x_k|)\one(|x_0|+\cdots+|x_k|\leq 1) \ dx \\
           & =\int_{[a,\infty)\times[b,\infty)} L_{k-1}(y+z)\one(y+z\leq 1)\,\nu_{l+1,a}(dy)\nu_{k-l,b}(dz).
\end{align*}
Here $\nu_{l+1,a}(dy)=S_l(a)\delta_a(dy)+\one(y>a)S_l(a,y)\,dy$, and similarly for $\nu_{k-l,b}$. The integration region is given by $y\geq a$, $z\geq b$ and $y+z\leq 1$. Expanding $\nu_{l+1,a}$ and $\nu_{k-l,b}$ gives four terms. Thus we have
\begin{align*}
    R(r,s) & = \int_{y=a}^{1-b} \int_{z=b}^{1-y} L_{k-1}(y+z) S_l(a,y) S_{k-l-1}(b,z) \ dzdy \\
           & \quad + S_{k-l-1}(b)\int_{y=a}^{1-b} L_{k-1}(y+b) S_l(a,y) \ dy \\
           & \quad + S_l(a)\int_{z=b}^{1-a} L_{k-1}(a+z)  S_{k-l-1}(b,z) \ dz \\
           & \quad + L_{k-1}(a+b) S_l(a) S_{k-l-1}(b).
\end{align*}
Observe that this is a polynomial in $r,s$, so let $\widetilde{R}$ be this polynomial, so that $R(r,s)=\widetilde{R}(r,s)$ for $1\geq r\geq |s|$. By symmetry, $R(r,s)=R(s,r)=R(-r,-s)$, so $R$ is a piecewise polynomial in $(-1,1)^2$. To prove (1), namely that $R$ is actually a polynomial, it suffices to show that $\widetilde{R}(r,s)=\widetilde{R}(s,r)=\widetilde{R}(-r,-s)$. Let
\begin{align*}
    T_1 & = \int_{y=a}^{1-b} \int_{z=b}^{1-y} L_{k-1}(y+z) S_l(a,y) S_{k-l-1}(b,z) \ dzdy, \\
    T_2 & = S_{k-l-1}(b)\int_{y=a}^{1-b} L_{k-1}(y+b) S_l(a,y) \ dy, \\
    T_3 & = S_l(a)\int_{z=b}^{1-a} L_{k-1}(a+z)  S_{k-l-1}(b,z) \ dz, \\
    T_4 & = L_{k-1}(a+b) S_l(a) S_{k-l-1}(b)
\end{align*}
be the four terms above, so that $\widetilde{R}(r,s)=T_1+T_2+T_3+T_4$. We will now show separately that $\widetilde{R}(r,s)=\widetilde{R}(-r,-s)$ and $\widetilde{R}(r,s)=\widetilde{R}(s,r)$.

\noindent\\
\textbf{Showing that $\widetilde{R}(r,s)=\widetilde{R}(-r,-s)$.} This is the easier of the two equalities. Note that $a=\frac{r+s}{2}$ and $b=\frac{r-s}{2}$, so under the transformation $(r,s)\mapsto (-r,-s)$, we have $(a,b)\mapsto (-a,-b)$. We will show that $T_1,T_2,T_3,T_4$ remain unchanged. For $T_4$, we have
\begin{align*}
    T_4(-r,-s) & = L_{k-1}(-a-b) S_l(-a) S_{k-l-1}(-b) \\
               & = (-1)^{k-1}L_{k-1}(a+b)(-1)^{l} S_l(a)(-1)^{k-l-1} S_{k-l-1}(b) \\
               & = T_4(r,s).
\end{align*}
For $T_3$, we have
\begin{align*}
    T_3(-r,-s) & = S_l(-a)\int_{z=-b}^{1+a} L_{k-1}(-a+z)  S_{k-l-1}(-b,z) \ dz \\
               & = (-1)^l S_l(a)\int_{z=-b}^{1+a} (-1)^{k-1}L_{k-1}(a-z)  (-1)^{k-l} S_{k-l-1}(b,-z) \ dz \\
               & = -S_l(a)\int_{z=-1-a}^{b} L_{k-1}(a+z) S_{k-l-1}(b,z) \ dz.
\end{align*}
Subtracting, we have
\begin{align*}
    T_3(r,s)-T_3(-r,-s) & = S_l(a)\int_{z=-1-a}^{1-a} L_{k-1}(a+z)  S_{k-l-1}(b,z) \ dz \\
                        & = S_l(a)\int_{z=-1}^{1} L_{k-1}(z) S_{k-l-1}(b,z-a) \ dz \\
                        & = 0,
\end{align*}
by the orthogonality of Legendre polynomials, since $S_{k-l-1}(b,z-a)$ is a polynomial in $z$ of degree at most $k-2$. Thus $T_3(-r,-s)=T_3(r,s)$. Similarly, one can show that $T_2(-r,-s)=T_2(r,s)$. Finally, for $T_1$, we first use the substitution $w=y+z$ to write
\begin{align*}
    T_1(r,s) & = \int_{w=a+b}^{1} \int_{y=a}^{w-b} L_{k-1}(w) S_l(a,y) S_{k-l-1}(b,w-y) \ dy dw.
\end{align*}
Then, we have
\begin{align*}
    T_1(-r,-s) & = \int_{w=-a-b}^{1} \int_{y=-a}^{w+b} L_{k-1}(w) S_l(-a,y) S_{k-l-1}(-b,w-y) \ dy dw \\
               & = \int_{w=-a-b}^{1} \int_{y=-a}^{w+b} L_{k-1}(-w)S_l(a,-y)S_{k-l-1}(b,-(w-y)) \ dy dw \\
               & = \int_{w=-1}^{a+b} \int_{y=w-b}^{a} L_{k-1}(w)S_l(a,y)S_{k-l-1}(b,w-y) \ dy dw \\
               & = -\int_{w=-1}^{a+b} \int_{y=a}^{w-b} L_{k-1}(w) S_l(a,y) S_{k-l-1}(b,w-y) \ dy dw.
\end{align*}
Subtracting, we have
\begin{align*}
    T_1(r,s)-T_1(-r,-s) & = \int_{w=-1}^{1} \int_{y=a}^{w-b} L_{k-1}(w) S_l(a,y) S_{k-l-1}(b,w-y) \ dy dw \\
                        & = \int_{w=-1}^{1} L_{k-1}(w) \int_{y=a}^{w-b} S_l(a,y) S_{k-l-1}(b,w-y) \ dy dw \\
                        & = 0,
\end{align*}
again by the orthogonality of Legendre polynomials, noting that $S_l(a,y) S_{k-l-1}(b,w-y)$ is a polynomial in $w,y$ of degree at most $k-3$, hence $\int_{a}^{w-b} S_l(a,y) S_{k-l-1}(b,w-y) dy$ is of degree at most $k-2$. Thus $T_1(-r,-s)=T_1(r,s)$, hence $\widetilde{R}(-r,-s)=\widetilde{R}(r,s)$.

\noindent\\
\textbf{Showing that $\widetilde{R}(s,r)=\widetilde{R}(r,s)$.} Under the transformation $(r,s)\mapsto (s,r)$, we have $(a,b)\mapsto (a,-b)$. Let $\Delta T_i(r,s)=T_i(s,r)-T_i(r,s)$ for $i=1,2,3,4$. Let us first rewrite $T_1$ as
\begin{align*}
    T_1(r,s) & = \int_{z=b}^{1-a} \int_{y=a}^{1-z} L_{k-1}(y+z) S_l(a,y) S_{k-l-1}(b,z) \ dydz.
\end{align*}
Then we have
\begin{align*}
    T_1(s,r) & = \int_{z=-b}^{1-a} \int_{y=a}^{1-z} L_{k-1}(y+z) S_l(a,y) S_{k-l-1}(-b,z) \ dydz \\
             & = \int_{z=-b}^{1-a} \int_{y=a}^{1-z} L_{k-1}(y+z) S_l(a,y) S_{k-l-1}(b,z) \ dydz,
\end{align*}
so the difference is
\begin{align*}
    \Delta T_1(r,s) & = \int_{z=-b}^{b} \int_{y=a}^{1-z} L_{k-1}(y+z) S_l(a,y) S_{k-l-1}(b,z) \ dydz \\
                    & = \int_{z=-b}^{b} \int_{y=a}^{1-z} L_{k-1}(-y-z) S_l(a,-y) S_{k-l-1}(b,-z) \ dydz \\
                    & = \int_{z=-b}^{b} \int_{y=-1-z}^{-a} L_{k-1}(y+z) S_l(a,y) S_{k-l-1}(b,z) \ dydz.
\end{align*}
Thus, this can also be written as
\begin{align*}
    \Delta T_1 & = \frac12\int_{z=-b}^{b} \left(\int_{y=-1-z}^{1-z} L_{k-1}(y+z) S_l(a,y) S_{k-l-1}(b,z) \ dy\right. \\
               & \quad \left.- \int_{y=-a}^{a} L_{k-1}(y+z) S_l(a,y) S_{k-l-1}(b,z) \ dy\right)dz \\
               & = \frac12\int_{z=-b}^{b} \left(\int_{y=-1}^{1} L_{k-1}(y) S_l(a,y-z) S_{k-l-1}(b,z) \ dy\right. \\
               & \quad \left.- \int_{y=-a}^{a} L_{k-1}(y+z) S_l(a,y) S_{k-l-1}(b,z) \ dy\right)dz \\
               & = -\frac12\int_{z=-b}^{b} \int_{y=-a}^{a} L_{k-1}(y+z) S_l(a,y) S_{k-l-1}(b,z) \ dy dz \\
               & = -\frac{1}{2 l!(k-l-1)!}\int_{z=-b}^{b} \int_{y=-a}^{a} L_{k-1}(y+z) a^{l-1} L_l'(y/a) b^{k-l-2}L_{k-l-1}'(z/b) \ dy dz \\
               & = -\frac{a^l b^{k-l-1}}{2 l!(k-l-1)!}\int_{z=-1}^{1} \int_{y=-1}^{1} L_{k-1}(ya+zb) L_l'(y)L_{k-l-1}'(z) \ dy dz. \\
\end{align*}
Letting $C=\frac{a^l b^{k-l-1}}{l!(k-l-1)!}$, we apply integration by parts to obtain
\begin{align*}
    \Delta T_1 & = -\frac12 C\int_{z=-1}^1 L_{k-l-1}'(z)\paren{[L_{k-1}(ya+zb) L_l(y) ]_{y=-1}^{1}-\int_{y=-1}^{1} L_{k-1}'(ya+zb) a L_l(y) \ dy} \ dz \\
               & = -\frac12 C\int_{z=-1}^1 L_{k-l-1}'(z)\paren{L_{k-1}(a+zb) - L_{k-1}(-a+zb)(-1)^l}dz \\
               & \quad +\frac12 Ca\int_{z=-1}^1\int_{y=-1}^{1} L_{k-l-1}'(z)L_{k-1}'(ya+zb) L_l(y) \ dy dz.
\end{align*}
Integrating by parts once more in the last term, we have
\begin{align*}
     & \int_{z=-1}^1\int_{y=-1}^{1} L_{k-l-1}'(z)L_{k-1}'(ya+zb) L_l(y) \ dy dz \\
     & = \int_{y=-1}^{1} L_l(y) \paren{[L_{k-l-1}(z)L_{k-1}'(ya+zb)]_{z=-1}^{1}}dy \\
     & \quad - \int_{z=-1}^{1}\int_{y=-1}^1 L_l(y)L_{k-l-1}(z) L_{k-1}''(ya+zb) b \ dydz \\
     & = \int_{y=-1}^{1} L_l(y) (L_{k-1}'(ya+b)-(-1)^{k-l-1}L_{k-1}'(ya-b))dy \\
     & \quad - \int_{z=-1}^{1}\int_{y=-1}^1 L_l(y)L_{k-l-1}(z) L_{k-1}''(ya+zb) b \ dydz. \\
\end{align*}
The second term is zero. Indeed, by expanding $L_{k-1}''(ya+zb)$ as a polynomial in $y,z$, each monomial is of the form $y^\alpha z^\beta$ with $\alpha+\beta\leq k-3$, so either $\alpha\leq l-1$ or $\beta\leq k-l-2$. If $\alpha\leq l-1$, then $\int_{-1}^{1} L_l(y) y^\alpha \ dy = 0$ by the orthogonality of Legendre polynomials. Similarly, if $\beta\leq k-l-2$, then $\int_{-1}^{1} L_{k-l-1}(z) z^\beta \ dz = 0$. Thus we have
\begin{align*}
    \Delta T_1 & = -\frac12 C\int_{z=-1}^1 L_{k-l-1}'(z)\paren{L_{k-1}(a+zb) - (-1)^l L_{k-1}(-a+zb)}dz \\
               & \quad + \frac12 Ca\int_{y=-1}^{1} L_l(y) (L_{k-1}'(ya+b)-(-1)^{k-l-1}L_{k-1}'(ya-b))dy.
\end{align*}
For $T_2$, we have
\begin{align*}
    T_2(s,r) & = S_{k-l-1}(-b)\int_{y=a}^{1+b} L_{k-1}(y-b) S_l(a,y) \ dy \\
             & = -S_{k-l-1}(b)\int_{y=-1-b}^{-a} L_{k-1}(y+b) S_l(a,y) \ dy,
\end{align*}
so the difference is
\begin{align*}
    \Delta T_2(r,s) & = -S_{k-l-1}(b)\paren{\int_{y=-1-b}^{1-b} L_{k-1}(y+b) S_l(a,y) \ dy - \int_{y=-a}^{a} L_{k-1}(y+b) S_l(a,y) \ dy} \\
                    & = -S_{k-l-1}(b)\paren{\int_{y=-1}^{1} L_{k-1}(y) S_l(a,y-b) \ dy - \int_{y=-a}^{a} L_{k-1}(y+b) S_l(a,y) \ dy} \\
                    & = S_{k-l-1}(b)\int_{y=-a}^{a} L_{k-1}(y+b) S_l(a,y) \ dy \\
                    & = \frac{C}{a}\int_{y=-a}^{a} L_{k-1}(y+b) L_l'(y/a) \ dy \\
                    & = C\int_{y=-1}^{1} L_{k-1}(ya+b) L_l'(y) \ dy.
\end{align*}
Integrating by parts gives
\begin{align*}
    \Delta T_2 & = C\paren{\sqb{L_{k-1}(ya+b) L_l(y)}_{y=-1}^{1} - \int_{y=-1}^{1} L_{k-1}'(ya+b) a L_l(y) \ dy} \\
               & = C(L_{k-1}(a+b) - (-1)^l L_{k-1}(-a+b)) - Ca\int_{y=-1}^{1} L_{k-1}'(ya+b) L_l(y) \ dy.
\end{align*}
For $T_3$, we have
\begin{align*}
    T_3(s,r) & = S_l(a)\int_{z=-b}^{1-a} L_{k-1}(a+z) S_{k-l-1}(-b,z) \ dz \\
             & = S_l(a)\int_{z=-b}^{1-a} L_{k-1}(a+z) S_{k-l-1}(b,z) \ dz,
\end{align*}
so the difference is
\begin{align*}
    \Delta T_3(r,s) & = S_l(a)\int_{z=-b}^{b} L_{k-1}(a+z) S_{k-l-1}(b,z) \ dz \\
                    & = C\int_{z=-1}^{1} L_{k-1}(a+zb) L_{k-l-1}'(z) \ dz.
\end{align*}
For $T_4$, we have
\begin{align*}
    T_4(s,r) & = L_{k-1}(a - b) S_l(a) S_{k-l-1}(-b) \\
             & = (-1)^{k-l-1} L_{k-1}(a - b) S_l(a) S_{k-l-1}(b),
\end{align*}
so the difference is
\begin{align*}
    \Delta T_4(r,s) & = \paren{(-1)^{k-l-1} L_{k-1}(a - b) - L_{k-1}(a + b)} S_l(a) S_{k-l-1}(b) \\
                    & = C\paren{(-1)^{k-l-1} L_{k-1}(a - b) - L_{k-1}(a + b)}.
\end{align*}

Summing up, we have
\begin{align*}
    \widetilde{R}(s,r)-\widetilde{R}(r,s)
     & = \Delta T_1 + \Delta T_2 + \Delta T_3 + \Delta T_4 \\
     & = -\frac12 C\int_{z=-1}^1 L_{k-l-1}'(z)\paren{L_{k-1}(a+zb) - (-1)^l L_{k-1}(-a+zb)}dz \\
     & \quad + \frac12 Ca\int_{y=-1}^{1} L_l(y) (L_{k-1}'(ya+b)-(-1)^{k-l-1}L_{k-1}'(ya-b))dy \\
     & \quad + C(L_{k-1}(a+b) - (-1)^l L_{k-1}(-a+b)) - Ca\int_{y=-1}^{1} L_{k-1}'(ya+b) L_l(y) \ dy \\
     & \quad + C\int_{z=-1}^{1} L_{k-1}(a+zb) L_{k-l-1}'(z) \ dz \\
     & \quad + C\paren{(-1)^{k-l-1} L_{k-1}(a - b) - L_{k-1}(a + b)} \\
     & = \frac12 C\int_{z=-1}^1 L_{k-l-1}'(z)\paren{L_{k-1}(a+zb) + (-1)^l L_{k-1}(-a+zb)}dz \\
     & \quad - \frac12 Ca\int_{y=-1}^{1} L_l(y) (L_{k-1}'(ya+b) + (-1)^{k-l-1}L_{k-1}'(ya-b))dy.
\end{align*}
Both integrals are zero. Indeed, for the first integral, we have
\[L_{k-1}(a+zb) + (-1)^l L_{k-1}(-a+zb) = L_{k-1}(a+zb) + (-1)^{k-l-1} L_{k-1}(a-zb),\]
which is even if $k-l-1$ is even and odd if $k-l-1$ is odd, while $L_{k-l-1}'(z)$ is odd if $k-l-1$ is even and even if $k-l-1$ is odd. Therefore their product is always an odd function, so the integral over $[-1,1]$ is zero. Similarly, for the second integral, we have
\[L_{k-1}'(ya+b) + (-1)^{k-l-1}L_{k-1}'(ya-b) = L_{k-1}'(ya+b) - (-1)^l L_{k-1}'(-ya+b),\]
which is odd if $l$ is even and even if $l$ is odd, while $L_l(y)$ is odd if $l$ is odd and even if $l$ is even. Therefore their product is always an odd function, so the integral over $[-1,1]$ is zero. Hence $\widetilde{R}(s,r)=\widetilde{R}(r,s)$. This proves (1) of the theorem.

It follows that $\deg_r \widetilde{R}=\deg_s \widetilde{R}$. To see that $\deg_s \widetilde{R}\leq k-1$, we show that $\deg_s T_i\leq k-1$ for $i=1,2,3,4$. For $T_4$, we have
\[T_4(r,s)=L_{k-1}(a+b) S_l(a) S_{k-l-1}(b)=L_{k-1}(r) S_l(a) S_{k-l-1}(b),\]
which is of degree at most $k-1$ in $s$. For $T_3$, we have
\[T_3(r,s)= S_l(a)\int_{z=b}^{1-a} L_{k-1}(a+z)  S_{k-l-1}(b,z) \ dz = S_l(a)\int_{z=r}^{1} L_{k-1}(z)  S_{k-l-1}(b,z-a) \ dz,\]
which is of degree at most $k-1$ in $s$. Similarly, for $T_2$, we have
\[T_2(r,s)= S_{k-l-1}(b)\int_{y=a}^{1-b} L_{k-1}(y+b) S_l(a,y) \ dy = S_{k-l-1}(b)\int_{y=r}^{1} L_{k-1}(y) S_l(a,y-b) \ dy,\]
which is of degree at most $k-1$ in $s$. Finally, for $T_1$, we have
\begin{align*}
    T_1(r,s) & = \int_{y=a}^{1-b} \int_{z=b}^{1-y} L_{k-1}(y+z) S_l(a,y) S_{k-l-1}(b,z) \ dzdy \\
             & = \int_{y=r}^{1} \int_{z=0}^{1-y} L_{k-1}(y+z) S_l(a,y-b) S_{k-l-1}(b,z+b) \ dzdy,
\end{align*}
which is of degree at most $k-1$ in $s$. This proves (2) of the theorem.

To show that $\int_{-1}^1 R(r,s) ds=2/k!$, we see that the left hand side is equal to $2\int_{V_r} Q(x) dx$, where $V_r$ is the affine subspace $x_0+\cdots+x_k=r$. Set $U(r)=\int_{V_r} Q(x) dx$. We want to show that $U(r)=1/k!$ for $0<r<1$. Just as in the computation of $R(r,s)$, we have
\begin{align*}
    U(r) & = \int_{V_r} Q(x) dx=\int_{V_r}L_{k-1}(|x_0|+\cdots+|x_k|)\one(|x_0|+\cdots+|x_k|\leq 1) \ dx \\
         & =\int_{y} L_{k-1}(y)\one(y\leq 1)\int_{x_0+\cdots+x_k=r}\delta(|x_0|+\cdots+|x_k|-y) \,dx\,dy \\
         & = S_k(r)L_{k-1}(r)+\int_{y=r}^1 L_{k-1}(y)S_k(r,y) \ dy \\
         & = \frac{1}{k!}\paren{r^k L_{k-1}(r)+r^{k-1}\int_{y=r}^1 L_{k-1}(y)L_k'(y/r) \ dy}.
\end{align*}
Observe that
\begin{align*}
    \int_r^1 L_{k-1}(y)L_k'(y/r) \ dy & = \int_r^1 L_{k-1}(-y)L_k'(-y/r) \ dy \\
                                      & = \int_{-1}^{-r} L_{k-1}(y)L_k'(y/r) \ dy,
\end{align*}
so we have
\begin{align*}
    \int_r^1 L_{k-1}(y)L_k'(y/r) \ dy & = \frac12\paren{\int_{-1}^{1} L_{k-1}(y)L_k'(y/r) \ dy - \int_{-r}^{r} L_{k-1}(y)L_k'(y/r) \ dy}.
\end{align*}
The leading term of $L_k(x)$ is $\frac{(2k)!}{2^k(k!)^2}x^k$, so the leading term of $L_k'(x)$ is $\frac{(2k-1)!}{2^{k-1}(k-1)!^2}x^{k-1}$. Thus $L_k'(y/r)$ can be written as $\frac{2k-1}{r^{k-1}} L_{k-1}(y)\ +$ lower degree terms in $y$. By the orthogonality of Legendre polynomials and the fact that $\int_{-1}^1 L_{k-1}(y)^2 \ dy=\frac{2}{2k-1}$, we have
\begin{align*}
    \int_{-1}^{1} L_{k-1}(y)L_k'(y/r) \ dy & = \frac{2k-1}{r^{k-1}}\int_{-1}^{1} L_{k-1}(y)^2 \ dy + 0 = \frac{2}{r^{k-1}}.
\end{align*}
For the second term, integrating by parts gives
\begin{align*}
    \int_{-r}^{r} L_{k-1}(y)L_k'(y/r) \ dy
     & = \sqb{L_{k-1}(y) r L_k(y/r)}_{y=-r}^{r} - r\int_{-r}^{r} L_{k-1}'(y) L_k(y/r) \ dy \\
     & = r L_{k-1}(r) L_k(1) - r L_{k-1}(-r) L_k(-1) - 0 \\
     & = 2r L_{k-1}(r).
\end{align*}
Together, we have
\begin{align*}
    U(r) & = \frac{1}{k!}\paren{r^k L_{k-1}(r) + \frac{r^{k-1}}{2}\paren{\frac{2}{r^{k-1}} - 2r L_{k-1}(r)}} = \frac{1}{k!}.
\end{align*}
This proves (3), and completes the proof of \Cref{thm:crosspoly}.

\subsection{Discrete approximation}

We have the following discrete approximation of \Cref{thm:crosspoly}.
\begin{lem} \label{lem:discretecrosspoly}
    Let $k$ be a fixed positive integer and $N$ be a large positive integer. Let
    \[\cP=\setcond{(x_0,\ldots,x_k)\in \frac12 \ZZ^{k+1}}{\sum |x_i|\leq N, \sum x_i\in \ZZ}.\]
    Then there exist a constant $C>0$ (depending only on $k$), a function $Q:\cP\to \RR$ and, for each subset $S\subsetneq [k]$, a polynomial $\widetilde{R}_S(r,s)$ such that the following holds. Given $S\subsetneq [k]$, for integers $r,s\in [-N,N]$, let $R_S(r,s)=\sum_{x\in\cP\cap V_{r,s}} Q(x)$, where $V_{r,s}$ is the affine subspace given by $x_0+\cdots+x_k=r$ and $x_0+\sum_{i\in S} x_i-\sum_{i\in [k]\setminus S}x_i=s$. Then:
    \begin{enumerate}
        \item For each $x\in\cP$, $|Q(x)|\leq CN^{-k}$.
        \item For each $r,s$, $|R_S(r,s)|\leq C/N$.
        \item $\sum_{s=-N}^N R_S(r,s)=1$.
        \item $\deg_r \widetilde{R}_S=\deg_s \widetilde{R}_S\leq k-1$.
        \item For all $r,s\in [-N,N]$, $|R_S(r,s)-\widetilde{R}_S(r,s)|\leq C/N^2$.
    \end{enumerate}
\end{lem}

\begin{proof}
    We will discretise the construction in \Cref{thm:crosspoly} to obtain $Q^*$, which will satisfy the required properties, except for (3), which we only have approximately. We can then renormalize $Q^*$ to obtain $Q$ satisfying all the properties.

    Let $Q_0:\RR^{k+1}\to \RR$ be given by $Q_0(x_0,\ldots,x_k)=L_{k-1}(|x_0|+\cdots+|x_k|)$ when $0\le |x_0|+\cdots+|x_k|\le 1$, and $Q_0(x)=0$ elsewhere. Then $Q_0$ is bounded, say $|Q_0(x)|\leq C$. Define $Q^*:\cP\to \RR$ and a polynomial $\widetilde{R}_S(r,s)$ by
    \[
        Q^*(x)=\frac{k!}{(2N)^k}Q_0(x/N),
        \qquad
        \widetilde{R}_S(r,s)=\frac{k!}{2N}R_{|S|}(r/N,s/N),
    \]
    where $R_l$ is as in \Cref{thm:crosspoly}. Then $\widetilde{R}_S$ satisfies (4). Without loss of generality, assume that $S=[l]$ for some $0\leq l\leq k-1$. Defining $R_S^*$ for $Q^*$, since $R_l$ is bounded, we have the second property $|R_S^*(r,s)|\leq C/N$. Then, we have
    \begin{align*}
        R_S^*(r,s)
         & = \frac{k!}{(2N)^k}
        \sum_{x\in\cP\cap V_{r,s}} Q_0(x/N) \\
         & = \frac{k!}{2N}\paren{\frac{1}{(2N)^{k-1}}\sum_{x\in\cP\cap V_{r,s}} Q_0(x/N)} \\
         & = \frac{k!}{2N}\paren{\int_{V_{r/N,s/N}} Q_0(y)\,dy+O(1/N)} \quad \text{(Riemann sum approximation)} \\
         & = \frac{k!}{2N}R_l(r/N,s/N)+O(1/N^2) \\
         & = \widetilde{R}_S(r,s)+O(1/N^2),
    \end{align*}
    giving (5). For the approximate version of (3), we have
    \begin{align*}
        \sum_{s=-N}^N R_S^*(r,s)
         & = \sum_{x_0+\cdots+x_k=r} Q^*(x) \\
         & = \frac{k!}{(2N)^k}\sum_{x_0+\cdots+x_k=r}Q_0(x/N) \\
         & = k!\int_{y_0+\cdots+y_k=r/N}Q_0(y)\,dy+O(1/N) \quad \text{(Riemann sum approximation)} \\
         & = 1+O(1/N).
    \end{align*}
    For each fixed level $r$, define
    \[Z_r:=\sum_{y: y_0+\cdots+y_k=r} Q^*(y)=1+O(1/N).
    \]
    Normalising with $Q(x)=Q^*(x)/Z_{x_0+\cdots+x_k}$, we get
    \[R_S(r,s)=(1+O(1/N))R_S^*(r,s),\]
    which gives (2)--(5). Finally, since $Z_r=1+O(1/N)$, the bound on $Q_0$ gives $|Q(x)|=O(N^{-k})$ uniformly on $\cP$, proving (1) after increasing $C$ if necessary.
\end{proof}

\begin{rmk}
    Based on numerical computations, we believe that an exact version of \Cref{lem:discretecrosspoly} holds. That is, there is a function $Q$ that is a polynomial in $|x_0|+\cdots+|x_k|$ so that $R_S=\widetilde{R}_S$. However, we do not have a proof, and the approximate version is sufficient for our application.
\end{rmk}

\section{Polynomial decomposition}

In this section, we prove the following result, which is \Cref{thm:decomp-rect} in the case that all $m_i=1$.

\begin{thm} \label{thm:decomp-linear}
    There exists $C>0$ such that the following holds. Let $F\in \ZZ[x]$ be a polynomial of degree at most $n$ such that $(1-x)^k\mid F$. Then $F$ can be written as a sum of at most $C\max(|F|,n)$ basic polynomials, each of degree at most $n$.
\end{thm}

The proof idea is to repeatedly apply focusing steps and cancelling steps, as discussed above.
In the first subsection we use the cross polytope distribution to prove the real filling result for basic polynomials.
We adapt this proof in the second subsection, combining it with a randomised rounding argument, to prove \Cref{thm:decomp-linear}.
The third subsection contains a sparse variant intended for future applications.
The final subsection contains a lower bound construction showing that there is no uniform filling result in the integer setting.

\subsection{Real filling}

In this subsection, we prove the following result, which is \Cref{thm:real} in the case that all $m_i=1$.

\begin{thm} \label{thm:real-basic}
    Fix $m=(1,\ldots,1) \in \ZZ_+^k$ and $n \in \ZZ_+$.
    Let $B$ be the set of polynomials in $\RR[x]$ of degree at most $n$.
    Let $\cB$ be the set of $m$-basic polynomials of degree at most $n$.
    Let $A = \RR^{\cB}$ and define $\pl: A \to B$ by $\pl(a) = \sum_{p \in \cB} a_p p$.
    Then there is a homomorphism $\delta:B \to A$ such that $\norm{\delta} \le C_{\ref{thm:real-basic}}$ for some $C_{\ref{thm:real-basic}}$ depending only on $k$, and for all $b \in \pl(A)$ we have
    $\pl \circ \delta(b)=b$.
\end{thm}

We first show that there is a homomorphism $\phi:B\to A$ that approximately inverts $\pl$, in the sense that $|b-\pl\circ \phi(b)|$ is small. Then to obtain $\delta$, we iteratively apply $\phi$ to the remainder and take the limit.

\begin{lem} \label{lem:real-approx}
    In the setting of \Cref{thm:real-basic}, there is a homomorphism $\phi:B \to A$ such that $\norm{\phi}\leq C$ for some $C$ depending only on $k$, and for all $F \in \pl(A)$, $\theta(F):=F-\pl \circ \phi(F)$ satisfies $|\theta(F)| \le C|F|/n$.
\end{lem}

\begin{proof}
    For small bounded $n$, we can take $\phi=0$. Henceforth, we assume that $n$ is sufficiently large in terms of $k$.

    We first note that the case of odd $n$ follows from the case of even $n$. Indeed, suppose $n=2n_0+1$. For a basic polynomial $p\in \cB$, write $[p]$ for the corresponding element in $A$. Define $\phi_0:B\to A$ by
    \[
        \phi_0(x^m)=\begin{cases}
            [x^{n-k}(x-1)^k] & \text{if }m=n,    \\
            0                & \text{otherwise,}
        \end{cases}
    \]
    and let $\phi_1:B\to A$ be the homomorphism given by the even case $2n_0$, with $\phi_1(x^{2n_0+1})=0$. Then we can define $\phi:B\to A$ by
    \[
        \phi(b)=\phi_0(b)+\phi_1(b-\pl\circ \phi_0(b)).
    \]
    Given $F\in \pl(A)$, the polynomial $F-\pl\circ \phi_0(F)\in \pl(A)$ has degree at most $2n_0$. Since $|\phi_0(b)|\leq |b|$ for all $b\in B$, we have
    \[
        |\phi(b)|\leq |\phi_0(b)|+O(|b-\pl\circ \phi_0(b)|)=O(|b|).
    \]
    Furthermore,
    \begin{align*}
        |F-\pl\circ \phi(F)| & = |F-\pl\circ \phi_0(F)-\pl\circ\phi_1(F-\pl\circ \phi_0(F))| \\
                             & = O(|F-\pl\circ \phi_0(F)|/n) = O(|F|/n).
    \end{align*}

    Henceforth, we assume that $n$ is even, which we will rewrite as $2n$. Let $F\in \pl(A)$. By multiplying $F$ by $x^{-n}$, we may work with a Laurent polynomial $F$ with exponents supported in $[-n,n]$. In this Laurent setting, we allow basic polynomials to have negative exponents, and ensure that we only use basic polynomials with exponents supported in $[-n,n]$. Afterwards, we shift everything back by multiplying by $x^n$ to obtain polynomials of degree at most $2n$. Thus, we redefine $B$ to be the set of polynomials in $\RR[x,x^{-1}]$ with exponents supported in $[-n,n]$ and $\cB$ to be the set of basic polynomials in $\RR[x,x^{-1}]$ with exponents supported in $[-n,n]$.

    We define $\phi:B\to A$ based on the weights of the discrete cross-polytope distribution as follows. Let $\cP,Q,\widetilde{R}_S,R_S$ be as in \Cref{lem:discretecrosspoly} with $N=n$. For a monomial $x^r$ with $-n\leq r\leq n$, let $\cP_r=\setcond{a\in \cP}{\sum a_i=r}$. Then
    \[
        \sum_{a\in \cP_r} Q(a)=1, \qquad \sum_{a\in \cP_r} |Q(a)|=O(1),
    \]
    where the second bound follows from $|Q(a)|=O(n^{-k})$ and $|\cP_r|=O(n^k)$. Set
    \[
        \phi(x^r)=\sum_{a\in \cP_r} Q(a)[x^r(1-x^{-2a_1})(1-x^{-2a_2})\cdots (1-x^{-2a_k})].
    \]

    First we check that each basic polynomial above has exponents supported in $[-n,n]$. By expanding these summands and looking at the exponents of their monomials, it suffices to show that for any $a\in \cP_r$ and $S\subseteq [k]$, we have $-n\leq r-2\sum_{j\in S} a_j\leq n$. Indeed, we have $a_0=r-\sum_{j=1}^k a_j$, so $\sum_{j=1}^k |a_j|+|r-\sum_{j=1}^k a_j|\leq n$. This implies that
    \begin{align*}
        \sum_{j\in S}a_j-\sum_{j\notin S}a_j-r+\sum_{j=1}^k a_j  & \leq n\quad \text{and} \\
        -\sum_{j\in S}a_j+\sum_{j\notin S}a_j+r-\sum_{j=1}^k a_j & \leq n.
    \end{align*}
    These yield $-n\leq r-2\sum_{j\in S} a_j\leq n$.

    Next, we check that $|\phi(F)|=O(|F|)$. Indeed, since $\sum_a |Q(a)|=O(1)$, we have $|\phi(x^r)|=O(1)$ for each $r$, and hence $|\phi(F)|=O(|F|)$.

    Lastly, let $\theta(F)=F-\pl\circ \phi(F)$. We will show that $|\theta(F)|=O(|F|/n)$. We first compute $\pl\circ \phi(x^r)$ for a monomial $x^r$. Write $a_S=\sum_{j\in S} a_j$ for $S\subseteq [k]$. Then we have
    \begin{align*}
        \pl\circ \phi(x^r)
         & = \sum_{a\in \cP_r} Q(a)x^r(1-x^{-2a_1})(1-x^{-2a_2})\cdots (1-x^{-2a_k}) \\
         & = x^r\sum_{a\in \cP_r} Q(a)\sum_{S\subseteq [k]} (-1)^{|S|}x^{-2a_S} \\
         & = x^r+\sum_{\emptyset\neq S\subseteq[k]} (-1)^{|S|}\sum_{a\in \cP_r} Q(a)x^{r-2a_S} \\
         & = x^r+\sum_{\emptyset\neq S\subseteq[k]} (-1)^{|S|}\sum_{-n\leq s\leq n} R_{[k]\setminus S}(r,s)x^{s} \\
         & = x^r+\sum_{\emptyset\neq S\subseteq[k]} (-1)^{|S|}\sum_{-n\leq s\leq n} (\widetilde{R}_{[k]\setminus S}(r,s)+O(n^{-2}))x^{s} \\
         & = x^r+\sum_{-n\leq s\leq n} (\widetilde{R}(r,s)+O(n^{-2}))x^{s}, \\
    \end{align*}
    where $\widetilde{R}(r,s)=\sum_{\emptyset\neq S\subseteq [k]} (-1)^{|S|}\widetilde{R}_{[k]\setminus S}(r,s)$ is a polynomial with $\deg_r \widetilde{R},\deg_s \widetilde{R}\leq k-1$.

    Write $F=\sum_{r=-n}^n c_r x^r$. Since $F\in \pl(A)$, we have $F(1)=F'(1)=\cdots=F^{(k-1)}(1)=0$, which implies that $\sum_r c_r r^j=0$ for all $j=0,\ldots,k-1$. Thus, we have
    \begin{align*}
        \theta(F) & = F-\pl\circ \phi(F) \\
                  & = -\sum_{-n\leq r\leq n} c_r\sum_{-n\leq s\leq n} (\widetilde{R}(r,s)+O(n^{-2}))x^{s} \\
                  & = -\sum_{-n\leq s\leq n}\paren{\sum_{-n\leq r\leq n} c_r(\widetilde{R}(r,s)+O(n^{-2}))}x^{s}.
    \end{align*}
    For each $s$, $\widetilde{R}(r,s)$ is a polynomial in $r$ of degree at most $k-1$, so $\sum_r c_r \widetilde{R}(r,s)=0$. Thus, we have
    \[
        \theta(F)=-\sum_{-n\leq s\leq n} \sum_{-n\leq r\leq n} c_r O(n^{-2})x^s=-\sum_{-n\leq s\leq n} O(|F|/n^2)x^s.
    \]
    Hence, $|\theta(F)|=O(|F|/n)$, as required.
\end{proof}

\begin{proof}[Proof of~\Cref{thm:real-basic}]
    We first prove the theorem for bounded $n$. For each $d=0,1,\ldots,n$, reduce the monomial $x^d$ modulo $(x-1)^k$ to obtain a polynomial $f_d$ of degree at most $k-1$. Then set $\delta(x^d)$ so that $\pl\circ \delta(x^d)=x^d-f_d$. Since $n$ is bounded, so is $\norm{\delta}$. Then for all $b\in B$, $b-\pl\circ \delta(b)$ is a polynomial of degree at most $k-1$. If $b\in \pl(A)$ is divisible by $(x-1)^k$, then so is $b-\pl\circ \delta(b)$, and hence $b-\pl\circ \delta(b)=0$.

    For larger $n$, let $\phi,\theta,C$ be as in \Cref{lem:real-approx}.
    First, let $T=\pl\phi:B\to \pl(A)$. By \Cref{lem:real-approx}, we have $\norm{1-T|_{\pl(A)}}\leq C/n<1$, so the series $U=\sum_{j=0}^\infty (1-T|_{\pl(A)})^j$ converges in operator norm to a bounded operator $\pl(A)\to \pl(A)$. Observe that $(1-T)U=U(1-T)=\sum_{j=1}^\infty (1-T)^j=U-1$, so $U=(T|_{\pl(A)})^{-1}$.

    Define the homomorphism $\delta$ by $\delta = \phi U^2 T$, which is bounded.
    Then for $b\in \pl(A)$, we have
    \[
        \pl\delta(b)=TU^2T(b)=b.
    \]
\end{proof}

\subsection{Integer decomposition}

Recall that a basic polynomial is of the form $\pm x^{a_0}(1-x^{a_1})\cdots (1-x^{a_k})$. The focusing procedure is the following.
\begin{lem}[Focusing step] \label{lem:focusing}
    Let $F\in \ZZ[x]$ be a polynomial of degree at most $n$. Then there exists $g$, a sum of at most $|F|$ basic polynomials each of degree at most $n$, such that $F-g$ is of degree at most $\ceil{\frac{k-1}{k}n}$ and $|F-g|\leq 2^k|F|$.
\end{lem}
\begin{proof}
    Set $m=\ceil{\frac{k-1}{k}n}$. For each positive integer $u\leq n$, let $g_u$ be the basic polynomial $x^u (1-x^{-v})^k$ where $v=\floor{u/k}$ is the largest integer for which $vk\leq u$. Besides the $x^u$ monomial of $g_u$, the next highest degree term is $x^{u-v}$, with $u-v=\ceil{\frac{k-1}{k}u}\leq m$. Write $F=\sum_u c_u x^u$, and let
    \[g=\sum_{m<u\leq n} c_u g_u.\]
    Then $g$ is a sum of at most $|F|$-many basic polynomials. $g$ eliminates all the monomials of $F$ of degree greater than $m$; thus $F-g$ is of degree at most $m$, and $|F-g|\leq 2^k|F|$.
\end{proof}

The cancelling step is roughly a random rounding of \Cref{thm:real-basic}.

\begin{lem}[Cancelling step] \label{lem:cancel}
    There exists a constant $C=C_{\ref{lem:cancel}}>0$ such that the following holds. Let $F\in \ZZ[x]$ be a polynomial of degree at most $n$ such that $(1-x)^k\mid F$ and $|F|=m$. Then there exists $G\in \ZZ[x]$, which can be written as a sum of at most $Cm$-many basic polynomials of degree at most $n$, such that $|F-G|\leq C\sqrt{mn}$.
\end{lem}

\begin{proof}
    Let $\cB$ be the set of basic polynomials of degree at most $n$. Let $\delta$ be the homomorphism from \Cref{thm:real-basic}, with $\norm{\delta}\leq C=C_{\ref{thm:real-basic}}$, where we may assume that $C>2$. Since $\pl\circ\delta(F)=F$, we may write
    \[
        F=\sum_{i=1}^t q_i p_i,
    \]
    where $p_i\in \cB$ are not necessarily distinct and $0<q_i\leq 1$, such that $\sum_i q_i\leq C|F|=Cm$. Pick independent $\varepsilon_i\in \{0,1\}$ with $\PP(\varepsilon_i=1)=q_i$ and take $G=\sum_i \varepsilon_i p_i$. We claim that with positive probability, at most $O(m)$ of the $\varepsilon_i$ are nonzero, and $|F-G|=O(\sqrt{mn})$, proving the lemma.

    Since $\EE[\sum_i \varepsilon_i]=\sum_i q_i\leq Cm$ and $\Var(\sum_i \varepsilon_i)=\sum_i q_i(1-q_i)\leq Cm$, we have by Chebyshev's inequality that
    \[
        \PP\paren{\sum_i \varepsilon_i\geq \EE\sqb{\sum_i \varepsilon_i}+Cm}\leq \frac{1}{Cm}.
    \]
    Thus with probability at least $1/2$, at most $2Cm$ of the $\varepsilon_i$ are nonzero.

    Since $\pl\circ \delta(F)=F$, we have
    \begin{equation} \label{eq:ex-f-g}
        \EE[F-G]=F-\sum_{i=1}^t q_i p_i=F-\pl\circ\delta(F)=0.
    \end{equation}
    Let $0\leq r\leq n$ and focus on the monomial $x^r$ of $F-G$. Let $\eta_i^{(r)}$ be the coefficient of $x^r$ in $p_i$. If $c_r$ is the coefficient of $x^r$ in $F$, then \eqref{eq:ex-f-g} says that $\EE[\sum_i \varepsilon_i \eta_i^{(r)}]=\sum_i q_i \eta_i^{(r)}=c_r$. The variance is
    \[
        \EE\sqb{\paren{\sum_i \varepsilon_i \eta_i^{(r)}-c_r}^2}=\sum_i q_i(1-q_i)(\eta_i^{(r)})^2\leq \sum_i q_i(\eta_i^{(r)})^2.
    \]
    Therefore,
    \begin{align*}
        \EE[|F-G|]^2 & = \EE\sqb{\sum_{r=0}^n \abs{c_r-\sum_{i=1}^t \varepsilon_i \eta_i^{(r)}}}^2 \\
                     & \leq (n+1)\sum_{r=0}^n \EE\sqb{\paren{c_r-\sum_{i=1}^t \varepsilon_i \eta_i^{(r)}}^2} \\
                     & \leq (n+1)\sum_{r=0}^n \sum_{i=1}^t q_i(\eta_i^{(r)})^2 \\
                     & = (n+1)\sum_{i=1}^t q_i\sum_{r=0}^n (\eta_i^{(r)})^2 \\
                     & \leq (n+1)\sum_{i=1}^t q_i\cdot O(1) \\
                     & =O(nm).
    \end{align*}
    By Markov's inequality, we have $|F-G|=O(\sqrt{nm})$ with probability, say, at least $2/3$. Therefore, with positive probability, both $|F-G|=O(\sqrt{nm})$ and at most $O(m)$ of the $\varepsilon_i$ are nonzero.
\end{proof}

We are now ready to prove \Cref{thm:decomp-linear} via focusing and cancelling.

\begin{proof}[Proof of \Cref{thm:decomp-linear}]
    Let $F\in \ZZ[x]$ be a polynomial of degree at most $n$ such that $(1-x)^k\mid F$. Define
    \[m=\max\paren{|F|,\frac{C_{\ref{lem:cancel}}^2 2^k (k+1)}{k}n}.\]
    Then our goal is to decompose $F$ into $O(m)$-many basic polynomials.

    We first consider a naive decomposition. Write $H=F/(1-x)^k\in\ZZ[x]$. By reading the monomials of $H$, we can write $F$ as a sum of $|H|$ basic polynomials, each of degree at most $n$. Using the long division algorithm, we have a simple bound of $|H|\leq 2^{O(n)}|F|$, enough to give $O(|F|)$ basic polynomials when $n=O(1)$.

    For larger $n$, we will apply the focusing and cancelling steps repeatedly, reducing the degree and the coefficient size of $F$ at each step, until we can apply the naive decomposition. By \Cref{lem:focusing}, there exists $g_1$, a sum of at most $m$-many basic polynomials, such that $F-g_1$ is of degree at most $\ceil{\frac{k-1}{k}n}\leq \frac{k}{k+1}n$ and $|F-g_1|\leq 2^k m$. By \Cref{lem:cancel}, there exists $g_2$, which can be written as a sum of at most $C_{\ref{lem:cancel}}\cdot 2^k m$-many basic polynomials, such that $F-g_1-g_2$ has degree at most $\frac{k}{k+1}n$ and
    \[
        |F-g_1-g_2|\leq C_{\ref{lem:cancel}} \sqrt{\frac{2^k kmn}{k+1}}\leq \frac{k}{k+1}m.
    \]

    Repeating the above process $O(\log n)$ times, we can reduce $F$ to a polynomial of degree at most $O(1)$ and size at most $O(m/n)$, using at most
    \[\sum_{i=0}^\infty C_{\ref{lem:cancel}}2^{k+1} \paren{\frac{k}{k+1}}^i m=O(m)\]
    many basic polynomials. Finally, we can apply the naive decomposition to write the remaining polynomial as a sum of $O(m/n)$-many basic polynomials.
\end{proof}

\subsection{Decomposition for sublinear-sized polynomials}

Using our decomposition result for polynomials of at least linear size, we can easily deduce a decomposition result for sublinear-sized polynomials,
which we record here for use in future applications.

\begin{cor} \label{cor:decomp-sublinear}
    Let $\varepsilon>0$. Then there exists $\delta>0$ such that the following holds. Let $F\in \ZZ[x]$ be a polynomial of degree at most $n$ such that $(1-x)^k\mid F$ and $|F|=O(n^{1-\varepsilon})$. Then $F$ can be written as a sum of $O(n^{1-\delta})$-many basic polynomials of degree at most $n$.
\end{cor}

\begin{proof}
    Assume that $n$ is sufficiently large. Suppose $F\in \ZZ[x]$ is of degree at most $n$ and $|F|\leq Cn^{1-\varepsilon}$ for some fixed $C>0$. By applying \Cref{lem:focusing} $m$ times ($m$ to be determined later), we obtain $g$, a sum of at most $2^{km}|F|$-many basic polynomials, such that $F-g$ is of degree at most $(\frac{k}{k+1})^m n$ and $|F-g|\leq 2^{km}|F|$. By picking $\delta>0$ sufficiently small, we can find $m$ such that $n^{\delta}<(\frac{k+1}{k})^m$ and $2^{km}<n^{\varepsilon-\delta}$. Then $\deg (F-g)\leq n^{1-\delta}$ and $|F-g|\leq Cn^{1-\delta}$. By \Cref{thm:decomp-linear}, $F-g$ is a sum of $O(n^{1-\delta})$ basic polynomials. Combining this with the decomposition of $g$, we see that $F$ is a sum of $O(n^{1-\delta})$ basic polynomials.
\end{proof}

\subsection{Lower bound} \label{subsec:lower}

In this section, we show that there is no uniform filling result in the integer setting.

\begin{prop}
    There exists $\epsilon>0$ such that the following holds. For $n$ sufficiently large, there exists a polynomial $F\in \ZZ[x]$ of degree at most $n$ such that $(1-x)^3\mid F$ and $|F|\leq 6$, but $F$ cannot be written as a sum of at most $\epsilon \log n$ basic polynomials.
\end{prop}

\begin{proof}
    Throughout the proof, we will ignore rounding issues for brevity. Consider the polynomial $F=(1-x^a)(1-x^b)-(1-x)(1-x^c)$ for some integers $a,b,c\geq 2$. Then $F(1)=F'(1)=0$, and $|F|\leq 6$. Furthermore, $F''(1)=0$ if and only if $ab=c$. Suppose $F$ can be written as a sum of at most $m$ basic polynomials for some $ab=c$ to be chosen later, with $m\leq \epsilon \log n$ for some sufficiently small $\epsilon>0$.
    \[
        (1-x^a)(1-x^b)-(1-x)(1-x^c)=\sum_{i=1}^m \varepsilon_i x^{d_i}(1-x^{a_i})(1-x^{b_i})(1-x^{c_i}).
    \]
    The left-hand side is the sum of the six monomials $x,-x^a,-x^b,x^{a+b},x^c,-x^{c+1}$, which we assume have distinct exponents. The right-hand side is the sum of $8m$ monomials, each with coefficient $\pm 1$ and an exponent that is a linear combination of the variables $a_i,b_i,c_i,d_i$. Thus for the above equation to hold, there exists a perfect matching between all $8m+6$ monomials, with each pair having the same exponent and suitable signs. Each pair asserts an affine linear equation in the variables $a,b,c,a_i,b_i,c_i,d_i$. For example, a pair within the $8m$ monomials on the right-hand side may give the equation
    \[
        d_i+a_i+b_i=d_j+a_j+c_j,
    \]
    whereas a pair between the left- and right-hand sides may give the equation
    \[
        d_i+b_i+c_i=c+1.
    \]

    Let $v=(a_1,b_1,c_1,d_1,a_2,\ldots,a_m,b_m,c_m,d_m)$ and write these $4m+3$ equations as
    \begin{equation} \label{eq:AB-linear}
        Av=B(a,b,c,1),
    \end{equation}
    where
    \begin{itemize}
        \item $A$ is a $(4m+3)\times 4m$ matrix with entries in $\{0,\pm 1\}$.
        \item Every row of $A$ has at most 8 non-zero entries.
        \item $B$ is a $(4m+3)\times 4$ matrix with entries in $\{0,\pm 1\}$.
    \end{itemize}
    Observe that an integer solution given by $v$ and $(a,b,c)$ to \eqref{eq:AB-linear} gives a decomposition of $F$ into at most $m$ basic polynomials; thus $ab=c$. This also holds over $\RR$, so in particular, not every $(a,b,c)\in \RR^3$ gives a real solution $v$ to \eqref{eq:AB-linear}.

    Let $r=\rank(A)$. Rearrange the rows and columns of $A$ (and the corresponding rows and columns of $B$ and $v$) so that the first $r$ rows and columns form an invertible $r\times r$ submatrix $C$. Let $\Delta=\det C\neq 0$. Every row $j>r$ of $A$ is thus a linear combination of the first $r$ rows, say
    \[
        A_j=\sum_{i=1}^r \alpha_{ji}A_i.
    \]
    By Cramer's rule,
    \[
        \alpha_{ji}=\frac{D_{ji}}{\Delta},
    \]
    where $D_{ji}$ is the determinant of $C$ with its $i$-th row replaced by the first $r$ entries of $A_j$. Since each row $A_j$ has at most 8 non-zero entries, Hadamard's inequality gives $|D_{ji}|\leq 8^{r/2}$. Similarly $|\Delta|\leq 8^{r/2}$.

    For $j>r$, the $j$-th row of \eqref{eq:AB-linear} gives
    \[
        B_j(a,b,c,1)=A_jv=\sum_{i=1}^r \alpha_{ji}A_iv=\sum_{i=1}^r \alpha_{ji}B_i(a,b,c,1).
    \]
    Thus
    \[
        \Delta B_j(a,b,c,1)=\Delta A_jv=\sum_{i=1}^r D_{ji}B_i(a,b,c,1),
    \]
    and these equations characterise exactly the $(a,b,c)\in \RR^3$ for which there exists a real solution $v$ to \eqref{eq:AB-linear}. Since not every $(a,b,c)\in \RR^3$ gives a real solution $v$, the above equation is non-trivial for some $j$. Say one of the equations is
    \[
        z_0+z_1a+z_2b+z_3c=0,
    \]
    where $z_i\in \ZZ$ are not all zero and $|z_i|\leq 8^{O(m)}$. For $\epsilon$ sufficiently small, we have $|z_i|< n^{1/5}$. Pick $(a,b,c)=(t,t^2,t^3)$ with $t=n^{1/4}$ large enough so that the above equation cannot be satisfied, and also $F$ has degree $c+1\leq n$, a contradiction.
\end{proof}

\section{Rectangular polynomial decomposition}

In this section, we prove \Cref{thm:decomp-rect} in full generality,
using our solution above to the special case that all $m_i=1$.

\subsection{Cyclotomic polynomials}

We begin by stating standard facts about cyclotomic polynomials. For more details, we direct the reader to \cite[Chapter 2]{Washington1997}.

For each positive integer $m$, the $m$-th cyclotomic polynomial $\Phi_m(x)$ is inductively defined by the formula
\[
    x^m-1=\prod_{d\mid m}\Phi_d.
\]
$\Phi_m$ is a monic, irreducible, integer polynomial of degree $\phi(m)$, where $\phi$ is the Euler totient function.

Let $\xi_m=e^{2\pi i/m}$, a root of the cyclotomic polynomial $\Phi_m$. Fix a basis $\ZZ[\xi_m]\cong \ZZ^{\phi(m)}$, including 1 as a basis element. For an element $x\in \ZZ[\xi_m]$, define the size of $x$, denoted $\norm{x}$, to be the $\ell^\infty$-norm of $x$ viewed as an element in $\ZZ^{\phi(m)}$.

This norm can be naturally extended to $\QQ[\xi_m]$. Observe that for $x,y\in \QQ[\xi_m]$, $\norm{xy}=O_m(\norm{x}\norm{y})$. If $0\neq y\in \ZZ[\xi_m]$ and $\norm{y}=O(1)$, then $y$ can only take boundedly many values, so $\norm{y^{-1}}=O(1)$.

We prove the following instance of \Cref{prob:general-ideal}, which allows us to transfer between ideals that differ by cyclotomic polynomials.

\begin{lem} \label{lem:add-cyclotomic}
    Fix a polynomial $R\in \ZZ[x]$ and integers $m\geq 1, a\geq 0$, such that $\Phi_m\nmid R$. There exists $C>0$ such that the following holds. Suppose $F\in \mathbb{Z}[x]$ is a polynomial of degree at most $n$ such that $\Phi_m^a R\mid F$. Then there exists $Q\in \ZZ[x]$ satisfying the following.
    \begin{enumerate}
        \item $\Phi_m^{a+1}R\mid F-Q$.
        \item $Q$ can be written as a sum of at most $C(|F|+\log n)$-many polynomials of degree at most $n$, each of which has one of the following forms:
        \begin{align}
             & \pm x^t \Phi_m^a R \quad\text{ for some } t\geq 0; \label{eq:add-cyclotomic-form1} \\
             & \pm x^t (x^{t_1m}-1)\cdots (x^{t_am}-1)R \quad \text{ for some } t,t_1,\ldots,t_a\geq 0. \label{eq:add-cyclotomic-form2}
        \end{align}
    \end{enumerate}
\end{lem}

\begin{proof}
    Observe that $\Phi_m^a\mid Q$ if and only if $Q^{(k)}(\xi_m)=0$ for $k=0,\ldots,a-1$. Thus, condition (1) is equivalent to having
    \begin{itemize}
        \item $R\mid Q$,
        \item $Q^{(k)}(\xi_m)=0$ for $k=0,\ldots,a-1$, and
        \item $Q^{(a)}(\xi_m)=F^{(a)}(\xi_m)$.
    \end{itemize}

    All implicit constants in this proof may depend on $R,m,a$, but not on $F$ or $n$.

    For bounded $n$, write $F=\Phi_m^a RH$ and take $Q=F$. Polynomial division by the fixed polynomial $\Phi_m^a R$ gives $|H|=O(|F|)$ many summands of the form \eqref{eq:add-cyclotomic-form1}, each of degree at most $n$. We may thus assume that $n$ is sufficiently large in terms of $R,m,a$.

    Let $P=\Phi_m^a R$. Then the following hold.
    \begin{itemize}
        \item $P^{(k)}(\xi_m)=0$ for $k=0,\ldots,a-1$.
        \item $P^{(a)}(\xi_m)=a!\Phi_m'(\xi_m)^a R(\xi_m)\neq 0$.
        \item $P^{(a)}(\xi_m)\mid F^{(a)}(\xi_m)$.
        \item $\norm{F^{(a)}(\xi_m)}=O(n^a |F|)$.
    \end{itemize}

    Let $q=F^{(a)}(\xi_m)/P^{(a)}(\xi_m)$. Then $q\in \ZZ[\xi_m]$ and $\norm{q}=O(n^a|F|)$.

    First suppose that $a=0$. Expressing $q$ in the basis of $\ZZ[\xi_m]$, choose $S\in\ZZ[x]$ with $\deg S<\phi(m)$, $S(\xi_m)=q$ and $|S|=O(\norm{q})=O(|F|)$. Set $Q=SR$. Then $R\mid Q$ and $Q(\xi_m)=F(\xi_m)$, so $\Phi_m R\mid F-Q$. Expanding $S$ writes $Q$ as a sum of $O(|F|)$-many polynomials of the form \eqref{eq:add-cyclotomic-form1}, each of bounded degree.

    Henceforth assume that $a\geq 1$. Write $q=q_0+m^ar$ where $\norm{q_0}\leq m^a$ and $\norm{r}=O(n^a|F|)$. By writing the coordinates of $r$ in base $2^a$ up to $b$ digits, we can express
    \[
        r=r_0+2^ar_1+\cdots+2^{ab} r_b,
    \]
    where $\norm{r_j}\leq 2^a$ for $j=0,1,\ldots,b-1$, and $\norm{r_b}=O(n^a|F|/2^{ab})$. Take $b$ to be the largest integer such that $2^b\leq \frac{n}{2ma}$. Then $2^b=\Theta(n)$, so $\norm{r_b}=O(|F|)$.

    By writing $q_0$ in terms of powers of $\xi_m$, we can pick a polynomial $S\in \ZZ[x]$ of degree at most $\phi(m)$ such that $S(\xi_m)=q_0$ and $|S|=O(1)$. Similarly, pick polynomials $S_j\in \ZZ[x]$ of degree at most $\phi(m)$ such that $S_j(\xi_m)=r_j$ for $j=0,1,\ldots,b$. Then $|S_j|=O(1)$ for $j=0,1,\ldots,b-1$ and $|S_b|=O(\norm{r_b})=O(|F|)$. Set
    \[
        Q=SP+\sum_{j=0}^b S_j x^a(x^{2^j m}-1)^a \Phi_m'^aR.
    \]
    The term $SP$ has size and degree $O(1)$ and can be written as a sum of $|S|$-many polynomials of the form \eqref{eq:add-cyclotomic-form1}. For $j=0,\ldots,b-1$, each term $S_j x^a(x^{2^j m}-1)^a \Phi_m'^aR$ has size $O(1)$ and degree at most $2^bma+O(1)\leq n$, and is a sum of $O(1)$-many polynomials of the form \eqref{eq:add-cyclotomic-form2}, obtained by expanding $S_jx^a\Phi_m'^a$ into monomials. There are $b=O(\log n)$ such terms. The final term $S_b x^a(x^{2^b m}-1)^a \Phi_m'^aR$ has size $O(|F|)$ and degree at most $2^bma+O(1)\leq n$, and the same expansion gives $O(|F|)$ summands of the form \eqref{eq:add-cyclotomic-form2}. Every summand also has degree at most $n$. Altogether, this gives $O(|F|+\log n)$ summands, as required.

    Note that $P\mid Q$ since $\Phi_m\mid x^{2^j m}-1$. We now verify that $Q^{(a)}(\xi_m)=F^{(a)}(\xi_m)$. Indeed, we have
    \begin{align*}
        Q^{(a)}(\xi_m) & = S(\xi_m)P^{(a)}(\xi_m)+\sum_{j=0}^b S_j(\xi_m)\paren{\left.\frac{d^a}{dx^a}\right|_{\xi_m}(x^{2^j m}-1)^a}\xi_m^a \Phi_m'(\xi_m)^a R(\xi_m) \\
                       & = S(\xi_m)a!\Phi_m'(\xi_m)^a R(\xi_m)+\sum_{j=0}^b S_j(\xi_m)a!(2^j m\xi_m^{2^j m})^a\Phi_m'(\xi_m)^a R(\xi_m) \\
                       & = a!\Phi_m'(\xi_m)^a R(\xi_m)\paren{q_0+\sum_{j=0}^b r_j 2^{ja} m^a} \\
                       & = a!\Phi_m'(\xi_m)^a R(\xi_m)q=F^{(a)}(\xi_m).
    \end{align*}
    Therefore, $\Phi_m^{a+1}R=\Phi_m P\mid F-Q$, as required.
\end{proof}

\subsection{Proof of \Cref{thm:decomp-rect}}

The proof goes by repeatedly applying \Cref{lem:add-cyclotomic} to transfer from the ideal $(\prod_i (1-x^{m_i}))$ to the ideal $((1-x^M)^k)$, where $M=\lcm(m_1,\ldots,m_k)$. Then we apply \Cref{thm:decomp-linear} to obtain the desired decomposition.

\begin{proof}[Proof of \Cref{thm:decomp-rect}]
    Let $M=\lcm(m_1,\ldots,m_k)$. We can write the ratio $(x^M-1)^k/\prod_i (x^{m_i}-1)$ as a product of cyclotomic polynomials, say $\Phi_{d_1}\cdots\Phi_{d_l}$. Applying \Cref{lem:add-cyclotomic} cumulatively for each $\Phi_{d_j}$, we obtain $Q_1,\ldots,Q_l\in\ZZ[x]$ such that for $j=1,\ldots,l$, we have
    \[
        \prod_{i=1}^k (x^{m_i}-1)\Phi_{d_1}\cdots\Phi_{d_j}\mid F-Q_1-\cdots-Q_j.
    \]
    Each $Q_j$ satisfies $|Q_j|=O(|F|+n)$ and $\deg Q_j\leq n$, and can be written as a sum of $O(|F|+n)$-many polynomials of the indicated forms. Setting $Q=Q_1+\cdots+Q_l$, we have $(x^M-1)^k\mid F-Q$, $|Q|=O(|F|+n)$ and $\deg Q\leq n$. We claim that $Q$ can be written as a sum of $O(|F|+n)$-many $(m_1,\ldots,m_k)$-basic polynomials, each of degree at most $n$. It suffices to establish this for each $Q_j$.

    Consider some $Q_j$ corresponding to $\Phi_{d_j}$. Let $P=\prod_{i=1}^k (x^{m_i}-1)\Phi_{d_1}\cdots \Phi_{d_{j-1}}$, and let $a$ be the largest exponent such that $\Phi_{d_j}^a\mid P$. From \Cref{lem:add-cyclotomic}, $Q_j$ can be written as a sum of $O(|F|+n)$-many polynomials of degree at most $n$, each of one of the following forms:
    \begin{itemize}
        \item $\pm x^t P$ for some $t\geq 0$;
        \item $\pm x^t (x^{t_1d_j}-1)\cdots (x^{t_ad_j}-1)P/\Phi_{d_j}^a$ for some $t,t_1,\ldots,t_a\geq 0$.
    \end{itemize}
    In the first case, expanding the fixed polynomial $P/\prod_i(x^{m_i}-1)$ into monomials gives $O(1)$-many $(m_1,\ldots,m_k)$-basic polynomials. In the second case, let $I=\setcond{i\in[k]}{d_j\mid m_i}$ and $h=|I|$. Since $h\leq a$, we may match each $i\in I$ with a distinct factor among the $a$ factors $x^{t_ud_j}-1$. For each matching corresponding to $m_i$ and a non-negative integer $t$, let $s$ be the largest integer such that $sm_i\leq td_j$. Then we can write
    \begin{align}
        \frac{(x^{td_j}-1)(x^{m_i}-1)}{\Phi_{d_j}}
         & = \frac{(x^{sm_i}-1+x^{td_j}-x^{sm_i})(x^{m_i}-1)}{\Phi_{d_j}} \nonumber \\
         & = \frac{x^{m_i}-1}{\Phi_{d_j}}(x^{sm_i}-1)+\frac{x^{td_j}-x^{sm_i}}{\Phi_{d_j}}(x^{m_i}-1). \label{eq:match-split}
    \end{align}
    Both $\frac{x^{m_i}-1}{\Phi_{d_j}}$ and $\frac{x^{td_j}-x^{sm_i}}{\Phi_{d_j}}$ are integer polynomials of size $O(1)$, where the size of the latter is due to $0\leq td_j-sm_i<m_i$.

    Now $P/\Phi_{d_j}^a$ can be written as the product of a fixed integer polynomial, the factors $x^{m_i}-1$ for $i\notin I$, and the factors $(x^{m_i}-1)/\Phi_{d_j}$ for $i\in I$. Rewrite each matched pair $(x^{t_ud_j}-1)(x^{m_i}-1)/\Phi_{d_j}$ using \eqref{eq:match-split}. This gives $O(1)$-many $(m_1,\ldots,m_k)$-basic polynomials for each original summand, without increasing its degree. Consequently, each $Q_j$, and hence $Q$, is a sum of $O(|F|+n)$-many such basic polynomials of degree at most $n$.

    Set $H=F-Q$, so that $(x^M-1)^k\mid H$, $\deg H\leq n$ and $|H|=O(|F|+n)$. Splitting the coefficients into powers of $x$ modulo $M$, write
    \[
        H(x)=\sum_{j=0}^{M-1}x^jH_j(x^M),
    \]
    where each $H_j$ is divisible by $(x-1)^k$. For each nonzero $H_j$, let $n_j=\floor{(n-j)/M}$. Then $\deg H_j\leq n_j$ and $|H_j|\leq |H|=O(|F|+n)$. By \Cref{thm:decomp-linear}, $H_j$ is a sum of $O(\max(|H_j|,n_j))=O(|F|+n)$ basic polynomials, each of degree at most $n_j$. Therefore, $x^jH_j(x^M)$ is a sum of $O(|F|+n)$ $(M,\ldots,M)$-basic polynomials of degree at most $j+Mn_j\leq n$. These are also $(m_1,\ldots,m_k)$-basic since $m_i\mid M$. Together with the decomposition of $Q$, this proves the theorem.
\end{proof}

As in \Cref{cor:decomp-sublinear}, we can also deduce an analogous decomposition result for sublinear-sized polynomials.

\begin{cor} \label{cor:decomp-rect-sublinear}
    Fix positive integers $m_1,\ldots,m_k$. For any $\varepsilon>0$, there exists $\delta>0$ such that the following holds. Let $F\in \ZZ[x]$ be a polynomial of degree at most $n$ such that $\prod_i (1-x^{m_i})\mid F$ and $|F|=O(n^{1-\varepsilon})$. Then $F$ can be written as a sum of $O(n^{1-\delta})$-many $(m_1,\ldots,m_k)$-basic polynomials. Furthermore, each such polynomial has degree at most $n$.
\end{cor}

\begin{proof}
    Set $M=\lcm(m_1,\ldots,m_k)$. The same focusing argument as in \Cref{cor:decomp-sublinear} works here with minor changes. Instead of focusing using polynomials of the form $x^u(1-x^{-v})^k$ in \Cref{lem:focusing}, we use polynomials of the form $x^u(1-x^{-vM})^k$. Then, instead of applying \Cref{thm:decomp-linear}, we apply \Cref{thm:decomp-rect}.
\end{proof}

\section{A two-variable lemma}

The goal of this section is to prove a lemma to be used in \cite{KL26-2}. This section is independent of the rest of the paper, and can be skipped by most readers. The lemma is included in this paper instead of \cite{KL26-2} since it is the only part of the proof requiring substantial algebra, and its proof is similar to the proof of \Cref{lem:add-cyclotomic}.

In this section, we will mostly work with Laurent polynomials in two variables: $\ZZ[x^{\pm},y^{\pm}]=\ZZ[x,x^{-1},y,y^{-1}]$. For a monomial $x^ay^b$, its degree is $|a|+|b|$. For a Laurent polynomial $F\in \ZZ[x^{\pm},y^{\pm}]$, its degree is the maximum degree of its monomials. The size of $F$, denoted $|F|$, is the sum of the absolute values of its coefficients.

\begin{prop} \label{prop:two-var}
    Fix a positive integer $r$ and vectors $u_i=(a_i,b_i)\in \ZZ^2$ with $a_i\neq 0$ for $i=1,\ldots,r$. Fix polynomials $A_i,B_i\in \ZZ[t]$ which are products of (not necessarily distinct) cyclotomic polynomials, with $A_i\mid B_i$ for all $i$, and let $Q\in\ZZ[t]$ be a product of distinct cyclotomic polynomials. In $\ZZ[x^{\pm},y^{\pm}]$, put
    \[
        p=\prod_{i=1}^r A_i(x^{a_i}y^{b_i}),\qquad p^*=\prod_{i=1}^r B_i(x^{a_i}y^{b_i}),\qquad q=Q(y).
    \]
    Then for sufficiently large $n$, if $F\in (p,q)$ is of degree at most $n$, there exists $G\in (p)$ such that $F-G\in (p^*,q)$ and $G$ is a sum of $O(|F|+\log n)$ Laurent polynomials in $(p)$, each of degree at most $n$ and size $O(1)$. All implicit constants may depend on the fixed data, but not on $F$ or $n$.
\end{prop}

In this section, let $Z=\ZZ[\zeta]$ be a fixed cyclotomic extension of $\ZZ$, where $\zeta=e^{2\pi i/m}$ is some root of unity. As before, fix an integral basis for $Z$, including 1 as a basis element, and let $\norm{\cdot}$ be the $\ell^\infty$ norm in this basis, which can be extended to $\QQ[\zeta]$. For $f=\sum_u c_u x^{u_1}y^{u_2}\in Z[x^{\pm},y^{\pm}]$, set $|f|=\sum_u\norm{c_u}$ and use the same definition of degree as above. Write $(h_1,\ldots,h_t)_Z$ for the ideal generated by $h_1,\ldots,h_t$ in $Z[x^{\pm},y^{\pm}]$.

\begin{lem}
    Let $f\in Z[x^{\pm},y^{\pm}]$. If $a,b\in Z$ are units and $f(a,b)=0$, then $f\in (x-a,y-b)_Z$.
\end{lem}
\begin{proof}
    For every integer $r$, $x^r-a^r$ is divisible by $x-a$ in $Z[x^{\pm}]$, and similarly for $y^r-b^r$. Replacing each monomial $x^r$ by $a^r$ and similarly for $y^r$, this shows that $f(x,y)-f(a,b)\in(x-a,y-b)_Z$.
\end{proof}

\begin{defn}
    For $f\in Z[x^{\pm},y^{\pm}]$ and units $a,b\in Z$, denote by $\ord_x f(a,b)$ the smallest non-negative integer $k$ such that $\partial_x^k f(a,b)\neq0$, if it exists. Otherwise, set $\ord_x f(a,b)=+\infty$.
\end{defn}

\begin{lem} \label{lem:ord+1}
    Let $f_1,\ldots,f_k\in Z[x^{\pm},y^{\pm}]$, and let $a,b\in Z$ be units. Suppose $f\in(f_1,\ldots,f_k)$, $\ord_x f_1(a,b)=\ell<\infty$, and $\ord_x f_j(a,b)>\ell$ for $j>1$. If $\ord_x f(a,b)>\ell$, then $f\in((x-a)f_1,(y-b)f_1,f_2,\ldots,f_k)$.
\end{lem}

\begin{proof}
    Write $f=g_1f_1+\cdots+g_k f_k$ for some $g_1,\ldots,g_k\in Z[x^{\pm},y^{\pm}]$. Then
    \begin{align*}
        0 = \partial_x^\ell f(a,b) = g_1(a,b)\partial_x^\ell f_1(a,b).
    \end{align*}
    Since $\partial_x^\ell f_1(a,b)\neq0$, we have $g_1(a,b)=0$. Thus, $g_1\in (x-a,y-b)$. Hence $f\in ((x-a)f_1, (y-b)f_1, f_2,\ldots,f_k)$.
\end{proof}

\begin{defn}
    A \emph{cyclotomic factor} over $Z$ is a Laurent polynomial of the form $x^\alpha y^\beta-\xi$, where $\alpha,\beta\in\ZZ$, and $\xi\in Z$ is a root of unity. A \emph{cyclotomic product} over $Z$ is a finite product of cyclotomic factors.
\end{defn}

For a cyclotomic factor $f(x,y)=x^\alpha y^\beta-\xi$ and non-zero $a,b\in\CC$ with $f(a,b)=0$, we have
\[\ord_x f(a,b)=\begin{cases}
        1      & \text{ if }\alpha\neq 0, \\
        \infty & \text{ if }\alpha=0.
    \end{cases}\]
Indeed, if $\alpha\neq0$, then $\partial_xf(a,b)=\alpha a^{-1}\xi\neq0$. If $\alpha=0$, then $f(x,b)$ is identically zero.

\begin{defn}
    Let $\cS$ be a class of Laurent polynomials over $\ZZ$ or a fixed cyclotomic extension $Z$. For a Laurent polynomial $f$, a \emph{$\cS$-decomposition of $f$ of size $m$} is an expression of $f$ of the form $f=f_1+\cdots+f_m$, where $f_1,\ldots,f_m\in\cS$, each of size $O(1)$. We say that the decomposition is \emph{$N$-bounded} if each $f_i$ has degree at most $N$.
\end{defn}

\begin{lem} \label{lem:match-partials}
    Let $Z$ be a fixed cyclotomic extension of $\ZZ$, and let $p\in Z[x^{\pm},y^{\pm}]$ be a fixed cyclotomic product over $Z$. Let $a,b\in Z$ be roots of unity and suppose $k=\ord_x p(a,b)<\infty$. Let $H\geq1$ and let $z\in Z$ be such that $\partial_x^k p(a,b)\mid z$ in $Z$ and $\norm{z}=O(Hn^k)$. Then, for sufficiently large $n$, there exists $g\in Z[x^{\pm},y^{\pm}]$ with an $n$-bounded $(p)_Z$-decomposition of size $O(H+\log n)$ such that $\partial_x^k g(a,b)=z$.
\end{lem}
\begin{proof}
    Write $p=\prod_{i=1}^{r}(x^{\alpha_i}y^{\beta_i}-\xi_i)$ and let $m_i$ be the order of $\xi_i$. Since $k=\ord_xp(a,b)<\infty$, exactly $k$ factors vanish at $(a,b)$, and each of these has $\alpha_i\neq0$, so we may assume that $a^{\alpha_i}b^{\beta_i}=\xi_i$ for $i\leq k$ and $a^{\alpha_i}b^{\beta_i}\neq\xi_i$ for $i>k$.

    First consider the case when $k=0$. Then $\norm{z}=O(H)$. Thus $\norm{z/p(a,b)}=O(H)$. Express $z/p(a,b)$ as $z_1+\cdots+z_s$ where each $\norm{z_t}=1$ and $s=O(H)$. Then take $g=\sum_{t=1}^s z_t p$, which satisfies $\partial_x^k g(a,b)=z$. This expression is an $O(1)$-bounded $(p)_Z$-decomposition of $g$ of size $O(H)$.

    Now consider the case when $k>0$. Let $p_1=\prod_{i=k+1}^{r}(x^{\alpha_i}y^{\beta_i}-\xi_i)$. Since $a^{\alpha_i}b^{\beta_i}=\xi_i$ for $i\leq k$, we have
    \[
        \partial_x^k p(a,b)=k!p_1(a,b)\prod_{i=1}^k(\alpha_i a^{-1}\xi_i)
        =k!p_1(a,b)a^{-k}\prod_{i=1}^k\alpha_i\xi_i.
    \]
    Since $\partial_x^k p(a,b)$ is fixed and non-zero, $\norm{\partial_x^k p(a,b)^{-1}}=O(1)$. Thus we have $\norm{z/\partial_x^k p(a,b)}=O(\norm{z}\norm{\partial_x^k p(a,b)^{-1}})=O(Hn^k)$. Let $M=\prod_{i=1}^k m_i\xi_i^{-1}$, and let $z_0\in Z$ have norm $O(1)$ such that
    \[
        M\mid\frac{z}{\partial_x^k p(a,b)}-z_0.
    \]
    Such a bounded choice of $z_0$ exists since $M$ is a fixed non-zero element of $Z$. Let
    \[
        z_1=\frac{z-z_0\partial_x^k p(a,b)}{M\partial_x^k p(a,b)}\in Z.
    \]
    Then $\norm{z_1}=O(Hn^k)$. Let $s$ be the largest integer such that
    \begin{equation} \label{eqn:sbound}
        2^s\sum_{i=1}^k (|\alpha_i|+|\beta_i|)m_i+\deg p_1+1<n.
    \end{equation}
    Then $2^s=\Theta(n)$. By writing the coordinates of $z_1$ in base $2^k$, we can express $z_1$ as $c_0+2^k c_1+\cdots+2^{ks}c_s$, where $c_i\in Z$, $\norm{c_i}=O(1)$ for $i=0,1,\ldots,s-1$, and $\norm{c_s}=O(H)$.

    Consider the polynomial
    \[g=z_0p+\sum_{t=0}^s c_t p_1\prod_{i=1}^k ((x^{\alpha_i}y^{\beta_i})^{2^t m_i}-1).\]
    Each factor $x^{\alpha_i}y^{\beta_i}-\xi_i$ for $i\leq k$ divides $(x^{\alpha_i}y^{\beta_i})^{2^tm_i}-1$, so every summand belongs to $(p)_Z$. To obtain a $(p)_Z$-decomposition of $g$ of size $O(H+\log n)$, it suffices to show that each summand has size $O(1)$ and there are $O(\log n)$ of them, except for the summand at $t=s$, which can be written as a sum of $O(H)$ polynomials of size $O(1)$.

    Indeed, $z_0p$ has size $O(1)$, and the polynomial
    \[c_t p_1\prod_{i=1}^k ((x^{\alpha_i}y^{\beta_i})^{2^t m_i}-1)\]
    has size $O(1)$ and belongs to $(p)_Z$ for $t=0,\ldots,s-1$, and can be written as a sum of $O(H)$ polynomials of size $O(1)$ in $(p)_Z$ for $t=s$, by writing $c_s$ as a sum of $O(H)$ elements of norm $1$. Furthermore, each polynomial in this decomposition has degree at most $n$ by \eqref{eqn:sbound}, so it is $n$-bounded. We now verify that $\partial_x^k g(a,b)=z$. We have
    \begin{align*}
        \partial_x^k g(a,b)
         & = z_0\partial_x^k p(a,b)+\sum_{t=0}^s c_t k!p_1(a,b) a^{-k}\prod_{i=1}^k 2^t m_i \alpha_i (a^{\alpha_i}b^{\beta_i})^{2^t m_i} \\
         & = z_0\partial_x^k p(a,b)+k!a^{-k}\sum_{t=0}^s c_t p_1(a,b) 2^{kt} \prod_{i=1}^k m_i \alpha_i \\
         & = z_0\partial_x^k p(a,b)+k!a^{-k}z_1 p_1(a,b) \prod_{i=1}^k m_i \alpha_i = z,
    \end{align*}
    as required.
\end{proof}

\begin{lem} \label{lem:add-cyclotomic-two-var}
    Let $Z$ be a fixed cyclotomic extension of $\ZZ$. Let $p\in Z[x^{\pm},y^{\pm}]$ be a fixed cyclotomic product over $Z$ such that each factor depends on $x$, and let $q=\prod_{j=1}^t(y-b_j)$ for distinct roots of unity $b_j\in Z$. Let $\varphi=x^\alpha y^\beta-\xi$ be a fixed cyclotomic factor over $Z$ with $\alpha\neq0$. Suppose $Z$ contains all roots $a$ of $\varphi(a,b_j)=0$ for every $j$. Then, for sufficiently large $n$, if $f\in(p,q)_Z$ has degree at most $n$, there exists $g\in Z[x^{\pm},y^{\pm}]$ with an $n$-bounded $(p)_Z$-decomposition of size $O(|f|+\log n)$ such that $f-g\in(\varphi p,q)_Z$.
\end{lem}

\begin{proof}
    If $t=0$, then $q=1$ and we may take $g=0$. Otherwise, set $H=|f|+\log n$. Replacing $\varphi$ by $x^{-\alpha}y^{-\beta}-\xi^{-1}$ if necessary, we may assume $\alpha>0$. Unless otherwise stated, all ideals in this proof are over $Z$. Let $Y=\{b_1,\ldots,b_t\}$, and for each $b\in Y$, let $X_b=\setcond{a\in Z}{a^\alpha b^\beta=\xi}$, so that $X_b$ consists of $\alpha$ distinct roots of unity.

    Enumerate $X_{b_j}=\{a_{1j},\ldots,a_{\alpha j}\}$. Instead of directly transferring from the ideal $(p,q)$ to $(\varphi p,q)$ in the lemma, we will do it in $\alpha t$ steps, where each step gets us closer to the desired ideal. This is done by constructing Laurent polynomials $g_{ij}\in Z[x^{\pm},y^{\pm}]$ for
    \[
        (i,j)=(1,1),(2,1),\ldots,(\alpha,1),(1,2),\ldots,(\alpha,t).
    \]
    For each $j$, set $f_{1j}=f_{\alpha,j-1}-g_{1j}$, where $f_{\alpha 0}=f$, and for $i\geq2$, set $f_{ij}=f_{i-1,j}-g_{ij}$. We will prove via induction that each $g_{ij}$ has an $n$-bounded $(p)_Z$-decomposition of size $O(H)$, each $f_{ij}$ has size $O(H)$ and degree at most $n$, and these Laurent polynomials belong to the ideal
    \begin{equation} \label{eqn:fij}
        f_{ij}\in \paren{\varphi p, p\prod_{l=1}^j(y-b_l), p\prod_{r=1}^i (x-a_{rj})\prod_{l=1}^{j-1} (y-b_l), q}_Z.
    \end{equation}
    In particular, since $\varphi(x,b_j)=b_j^\beta\prod_{r=1}^\alpha(x-a_{rj})$ and $b_j^\beta$ is a unit, we have
    \[
        ((x-a_{1j})\cdots(x-a_{\alpha j}),y-b_j)=(\varphi,y-b_j).
    \]
    Thus the case $i=\alpha$ simplifies to
    \begin{equation} \label{eqn:faj}
        f_{\alpha j}\in \paren{\varphi p, p\prod_{l=1}^j(y-b_l), q}_Z.
    \end{equation}

    The induction begins with $f_{\alpha 0}=f\in(p,q)_Z$, satisfying \eqref{eqn:faj}. For $i=1$, $j\geq1$, the induction hypothesis gives
    \[
        f_{\alpha,j-1}\in (\varphi p, (y-b_1)\cdots (y-b_{j-1})p, q)_Z.
    \]
    Let $p^*=(y-b_1)\cdots (y-b_{j-1})p$, and let
    \[
        k=\ord_x p(a_{1j},b_j)=\ord_x p^*(a_{1j},b_j).
    \]
    Write
    \[
        f_{\alpha,j-1}=A\varphi p+Bp^*+Cq,
    \]
    for some $A,B,C\in Z[x^{\pm},y^{\pm}]$. Since $q(b_j)=0$, $q$ is independent of $x$, and $\ord_x(\varphi p)(a_{1j},b_j)=k+1$, we get
    \[
        \partial_x^k f_{\alpha,j-1}(a_{1j},b_j)=B(a_{1j},b_j)\partial_x^k p^*(a_{1j},b_j).
    \]
    Thus,
    \[
        \partial_x^k p^*(a_{1j},b_j)\mid \partial_x^k f_{\alpha,j-1}(a_{1j},b_j).
    \]
    Since $|f_{\alpha,j-1}|=O(H)$ and $\deg f_{\alpha,j-1}\leq n$, we have
    \[
        \norm{\partial_x^k f_{\alpha,j-1}(a_{1j},b_j)}=O(Hn^k).
    \]
    By \Cref{lem:match-partials}, there exists $g_{1j}\in Z[x^{\pm},y^{\pm}]$ with an $n$-bounded $(p^*)_Z$-decomposition of size $O(H+\log n)=O(H)$ such that
    \[
        \partial_x^k g_{1j}(a_{1j},b_j)=\partial_x^k f_{\alpha,j-1}(a_{1j},b_j).
    \]
    Since $(p^*)_Z\subseteq(p)_Z$, this is also a $(p)_Z$-decomposition. Thus $|f_{1j}|=O(H)$ and $\deg f_{1j}\leq n$.

    Since $\ord_x f_{\alpha,j-1}(a_{1j},b_j)\geq k$, it follows that for $f_{1j}=f_{\alpha,j-1}-g_{1j}$, we have $\ord_x f_{1j}(a_{1j},b_j)\geq k+1$. We also have $f_{1j}\in (\varphi p, p^*, q)$, so by \Cref{lem:ord+1},
    \[
        f_{1j}\in (\varphi p, (x-a_{1j})p^*, (y-b_{j})p^*, q),
    \]
    establishing the induction step in this case.

    For $i>1,j\geq 1$, by our induction hypothesis, we have
    \[
        f_{i-1,j}\in (\varphi p, (y-b_1)\cdots (y-b_j)p, (x-a_{1j})\cdots (x-a_{i-1,j})(y-b_1)\cdots (y-b_{j-1})p, q).
    \]
    Let $p^*=(x-a_{1j})\cdots (x-a_{i-1,j})(y-b_1)\cdots (y-b_{j-1})p$. As above, for $k=\ord_x p(a_{ij},b_j)=\ord_x p^*(a_{ij},b_j)$, we can find $g_{ij}$ with an $n$-bounded $(p^*)_Z$-decomposition of size $O(H)$ such that $\partial_x^k g_{ij}(a_{ij},b_j)=\partial_x^k f_{i-1,j}(a_{ij},b_j)$. This is also a $(p)_Z$-decomposition. Then $|f_{ij}|=O(H)$ and $\deg f_{ij}\leq n$, and for $f_{ij}=f_{i-1,j}-g_{ij}$, we have $\ord_x f_{ij}(a_{ij},b_j)\geq k+1$. We also have $f_{ij}\in (\varphi p, (y-b_1)\cdots (y-b_j)p, p^*, q)$, so by \Cref{lem:ord+1}, we have
    \begin{align*}
        f_{ij} & \in (\varphi p, (y-b_1)\cdots (y-b_j)p, (x-a_{i,j})p^*, (y-b_j)p^*, q) \\
               & = (\varphi p, (y-b_1)\cdots (y-b_j)p, (x-a_{1j})\cdots (x-a_{i,j})(y-b_1)\cdots (y-b_{j-1})p, q).
    \end{align*}
    This completes the induction.

    Now set $g=\sum_{i,j}g_{ij}$. There are only $O(1)$ pairs $(i,j)$, so combining their decompositions gives an $n$-bounded $(p)_Z$-decomposition of $g$ of size $O(H)=O(|f|+\log n)$, and
    \[
        f-g=f_{\alpha t}\in\paren{\varphi p,\prod_{l=1}^t(y-b_l)p,q}_Z=(\varphi p,q)_Z,
    \]
    as required.
\end{proof}

We will also use the following observation. Let $\Psi:Z\to\ZZ$ be the $\ZZ$-linear map projecting onto the coefficient of $1$ in the chosen integral basis of $Z$. It fixes $\ZZ$ and satisfies $|\Psi(z)|\leq\norm{z}$. Extend $\Psi$ coefficient-wise to a $\ZZ$-linear map $\Psi:Z[x^{\pm},y^{\pm}]\to\ZZ[x^{\pm},y^{\pm}]$. Then this map is in fact $\ZZ[x^{\pm},y^{\pm}]$-linear.

\begin{lem} \label{lem:descent}
    Let $h_1,\ldots,h_k$ be polynomials in $\ZZ[x^{\pm},y^{\pm}]$. If $h\in (h_1,\ldots,h_k)_Z$, then
    \[
        \Psi(h)\in (h_1,\ldots,h_k)_\ZZ.
    \]
\end{lem}
\begin{proof}
    Since $h\in (h_1,\ldots,h_k)_Z$, we can write $h=\sum_{i=1}^k w_i h_i$ for some $w_i\in Z[x^{\pm},y^{\pm}]$. Then
    \[
        \Psi(h)=\sum_{i=1}^k \Psi(w_i h_i)=\sum_{i=1}^k \Psi(w_i) h_i\in (h_1,\ldots,h_k)_\ZZ. \qedhere
    \]
\end{proof}

In particular, if $p$ has integer coefficients, applying $\Psi$ to each summand of an $N$-bounded $(p)_Z$-decomposition gives an $N$-bounded $(p)_{\ZZ}$-decomposition.

\begin{proof}[Proof of \Cref{prop:two-var}]
    Choose a fixed cyclotomic extension $Z=\ZZ[\zeta]$ containing all roots needed to split the polynomials $A_i$, $B_i$ and $Q$, and all roots $a$ of $a^{a_i}b^{b_i}=\xi$ whenever $Q(b)=0$ and $B_i(\xi)=0$. These are finitely many roots of unity, since $a_i\neq0$. Over $Z$, both $p$ and $p^*$ are cyclotomic products whose factors depend on $x$, and we can write
    \[
        p^*=p\varphi_1\cdots\varphi_s,\qquad q=\prod_{j=1}^t(y-b_j),
    \]
    where each $\varphi_l$ is a cyclotomic factor depending on $x$ and the $b_j$ are distinct roots of unity.

    Set $p_0=p$ and $p_l=p\varphi_1\cdots\varphi_l$ for $l=1,\ldots,s$. Starting with $F_0=F$, apply \Cref{lem:add-cyclotomic-two-var} successively to obtain $g_l\in(p_{l-1})_Z$ and $F_l=F_{l-1}-g_l\in(p_l,q)_Z$. Each $g_l$ is a sum of $O(|F_{l-1}|+\log n)$ polynomials in $(p_{l-1})_Z$, each of size $O(1)$ and degree at most $n$. Since $s$ is fixed, induction gives
    \[
        |F_l|=O(|F|+\log n),\qquad \deg F_l\leq n.
    \]
    Thus $\widetilde G=\sum_{l=1}^s g_l$ is a sum of $O(|F|+\log n)$ polynomials in $(p)_Z$, each of size $O(1)$ and degree at most $n$, and
    \[
        F-\widetilde G=F_s\in(p^*,q)_Z.
    \]
    Finally, set $G=\Psi(\widetilde G)$. Since $F,p,p^*,q$ have integer coefficients and $\Psi$ fixes $F$, \Cref{lem:descent} gives $G\in(p)_\ZZ$ and $F-G\in(p^*,q)_\ZZ$.
\end{proof}

\section{Approximate counting} \label{sec:count}

We conclude our paper by proving our approximate counting result. This requires the following strengthening of \Cref{thm:decomp-rect}. We say that a monomial $x^r$ \emph{occurs} in an $m$-basic polynomial $x^a(1-x^{m_1a_1})\cdots(1-x^{m_ka_k})$ if $r=a+a_S$ for some subset $S \subseteq [k]$, where $a_S = \sum_{i \in S} m_ia_i$. A monomial may occur multiple times in the same polynomial, even if its coefficient is zero in the polynomial (since $x^r$ and $-x^r$ may both occur). Define the \emph{diameter} of a polynomial to be $b-a$, where $a$ and $b$ are the smallest and largest exponents of the monomials appearing in the polynomial. Call an $m$-basic polynomial \emph{simple} if all of its $2^k$ monomials have distinct exponents.

\begin{lem} \label{lem:decomp-bounded}
    Fix $m=(m_1,\ldots,m_k) \in \ZZ_+^k$.
    For any $\varepsilon>0$, there exists $C>0$ such that the following holds for sufficiently large $n$.
    Let $F\in \ZZ[x]$ be a polynomial of degree at most $n$ such that $(1-x^{m_1})\cdots (1-x^{m_k})\mid F$.
    Then $F$ can be written as a sum of at most $C\max(|F|,n)$ polynomials,
    each of which is $m$-basic, simple, has diameter at most $\varepsilon n$ and of degree at most $n$. Furthermore, each monomial occurs at most $C\norm{F}_\infty$ times in the decomposition.
\end{lem}

\begin{proof}
    We may assume that $\varepsilon\leq 1/4$. The case $F=0$ is trivial, so assume $\norm{F}_\infty>0$. By \Cref{thm:decomp-rect}, we have a decomposition of $F$ into $O(\max(|F|,n))$ $m$-basic polynomials of degree at most $n$.

    We first modify the decomposition so that each $m$-basic polynomial
    \[
        p=\pm x^a(1-x^{m_1a_1})\cdots (1-x^{m_ka_k})
    \]
    has $m_ia_i\leq\varepsilon n/(2k)$ for all $i$. If any of the $m_ia_i$ is larger than $\varepsilon n/(2k)$, using the identity
    \[
        1-x^{m_i(b+c)}=1-x^{m_ib}+x^{m_ib}(1-x^{m_ic}),
    \]
    we can break up the factor $1-x^{m_ia_i}$ into a sum of $O(1)$ summands of the form $x^d(1-x^{m_ib})$ with $1\leq m_ib\leq\varepsilon n/(2k)$. Repeating this for all $i$ increases our number of $m$-basic polynomials by only a constant factor, but now all the factors have $m_ia_i\leq\varepsilon n/(2k)$. Thus each summand has diameter at most $\varepsilon n/2$.

    Next, we make all the $m$-basic polynomials simple, while keeping their diameters at most $\varepsilon n$. Since $p$ has diameter at most $\varepsilon n/2\leq n/8$, its exponents are supported in either $[0,5n/8]$ or $[3n/8,n]$. Assume the former, as the latter follows by reflection. Choose integers $b_i$ with $m_ib_i\in(\varepsilon n/(2k),\varepsilon n/k)$. By replacing each factor $1-x^{m_ia_i}$ with $1-x^{m_ib_i}-x^{m_ia_i}(1-x^{m_i(b_i-a_i)})$, we can write $p$ as the sum of
    \[
        \pm(-1)^{|S|}x^{a+a_S}\prod_i(1-x^{m_i(b_i-1_{i\in S}a_i)}),
        \qquad S\subseteq[k].
    \]
    Each new polynomial has degree at most $n$ and diameter at most $\sum_i m_ib_i<\varepsilon n$. A collision between two monomials within any one of these polynomials imposes a non-trivial affine equation on the $b_i$. There are only $O(1)$ such equations, each excluding $O(n^{k-1})$ of the $\Theta(n^k)$ choices. Thus we may choose the $b_i$ so that all the resulting polynomials are simple. This increases the number of summands by only a constant factor.

    We describe a procedure to cancel a positive and a negative occurrence of $x^r$ in two distinct simple summands $p,q$. The replacement uses $O(1)$-many $m$-basic polynomials, removes the two selected occurrences, and adds at most $O(1)$ occurrences at each newly introduced exponent.

    Suppose $p$ contains $+x^r$ and $q$ contains $-x^r$. Write
    \begin{align*}
        p & =x^r(1-x^{m_1a_1})\cdots (1-x^{m_ka_k}),  & a_1,\ldots,a_k\in \ZZ, \\
        q & =-x^r(1-x^{m_1b_1})\cdots (1-x^{m_kb_k}), & b_1,\ldots,b_k\in \ZZ.
    \end{align*}
    Since $p,q$ have diameter at most $\varepsilon n$ and have $x^r$ in common, together they have diameter at most $2\varepsilon n\leq n/2$. Thus, both have exponents supported in either $[0,3n/4]$ or $[n/4,n]$. Assume the former, as the latter follows by reflection.
    Let $c_1,\ldots,c_k$ be integers with $0<m_ic_i<\varepsilon n/k$. Consider the $m$-basic polynomials
    \begin{align*}
        p_S & = (-1)^{|S|}x^{r+c_S}(1-x^{m_1a_1})\cdots (1-x^{m_ka_k}), \\
        q_S & = (-1)^{|S|+1}x^{r+c_S}(1-x^{m_1b_1})\cdots (1-x^{m_kb_k}), \\
        r_S & = (-1)^{|S|}x^{r+a_S}(1-x^{m_1c_1})\cdots (1-x^{m_kc_k}), \\
        s_S & = (-1)^{|S|+1}x^{r+b_S}(1-x^{m_1c_1})\cdots (1-x^{m_kc_k}).
    \end{align*}
    Then $p=p_\emptyset$, $q=q_\emptyset$, and $r_\emptyset=-s_\emptyset$. We have the identity
    \[
        \sum_{S\subseteq [k]} p_S=x^r(1-x^{m_1a_1})\cdots (1-x^{m_ka_k})(1-x^{m_1c_1})\cdots (1-x^{m_kc_k})=\sum_{S\subseteq [k]} r_S,
    \]
    and similarly $\sum_{S\subseteq [k]} q_S=\sum_{S\subseteq [k]} s_S$. Thus, replacing $p+q$ in the decomposition with
    \[
        \sum_{\emptyset\neq S\subseteq [k]} (-p_S+r_S-q_S+s_S)
    \]
    does not change the sum. There are two kinds of monomials (which we will refer to by their exponents) occuring in this replacement. First are the original monomials $r+a_S$ and $r+b_S$ from $p$ and $q$, now occuring in $r_S$ and $s_S$, respectively, with the exception of the two selected occurrences at $r$ corresponding to $S=\emptyset$.
    Monomials of the second kind are of the form $r+c_S+a_{S'}$ or $r+c_S+b_{S'}$ for some $S,S'\subseteq [k]$ with $S\neq \emptyset$, and they come in cancelling pairs, possibly with several pairs at the same exponent. For instance, $r+c_S+a_{S'}$ occurs in $p_S$ and $r_{S'}$ for $S'\neq\emptyset$, and $r+c_S$ occurs in $p_S$ and $q_S$.
    Each monomial of the second kind can occur at most $M:=2^{k+2}$ times in the replacement, and is a non-constant affine function of the $c_i$.
    Let $W\subseteq\{0,\ldots,n\}$ be a set of at most $\eta n$ exponents, where $\eta>0$ is sufficiently small. We require the monomials of the second kind to avoid $W$ and the supports of $p,q$.
    We also require that the $c_S$ be distinct for distinct subsets $S\subseteq[k]$, which adds only $O(1)$ non-trivial affine equations for the $c_i$ to avoid.
    Thus for sufficiently small $\eta$ and sufficiently large $n$, some choice of $c_i$ avoids all such collisions. Indeed, there are $\Theta((\varepsilon n)^k)$ possible choices for the $c_i$, and each affine equation forbids $O((\varepsilon n)^{k-1})$ of them. There are $O(|W|+1)=O(\eta n+1)$ equations to avoid, so choosing $\eta$ sufficiently small in terms of $\varepsilon$ leaves a valid choice of $c_i$. The polynomials $p_S,q_S$ are simple since they are translates of $p,q$, and $r_S,s_S$ are simple since the $c_S$ are distinct. Since $0<m_ic_i<\varepsilon n/k$, each new $m$-basic polynomial has degree at most $n$ and diameter at most $\varepsilon n$.

    To finish the proof, set $D=\norm{F}_\infty\geq 1$ and $B=LD$, where $L\geq 2M$ is a sufficiently large integer. Suppose each monomial $x^r$ occurs $t_r$ times in the decomposition. Let
    \[
        T=\sum_r \max(0,t_r-B)
    \]
    be the number of excess occurences, which is at most $\sum_r t_r=O(|F|+n)$. If $T=0$, then we are done. Otherwise, we will use steps that strictly decrease $T$ and introduce at most $O(1)$ summands. After at most $O(|F|+n)$ steps, we have $T=0$ and a decomposition of size $O(|F|+n)$. So we may assume that throughout this process, $\sum_r t_r\leq C_1(|F|+n)$ for some $C_1>0$.

    Since $T>0$, some monomial $x^r$ occurs more than $B$ times. We can find a pair of $m$-basic polynomials $p,q$ such that $p$ contains $x^r$ and $q$ contains $-x^r$. Define
    \[
        W=\setcond{s}{t_s>B-M}.
    \]
    We have $|W|\leq C_1(|F|+n)/(B-M)\leq 10C_1Dn/B\leq\eta n$ if $L$ is sufficiently large. Perform the replacement above avoiding $W$. The count at $r$ decreases by two, and all other occurrences from $p,q$ are retained. Each newly introduced exponent gains at most $M$ occurrences without exceeding the threshold $B$. Thus $T$ decreases by at least one.

    Repeating this process eventually reduces $T$ to zero, resulting in a decomposition of size $O(\max(|F|,n))$ into simple $m$-basic polynomials of degree at most $n$ and diameter at most $\varepsilon n$, in which each monomial occurs at most $B=O(\norm{F}_\infty)$ times.
\end{proof}

\begin{proof}[Proof of Theorem \ref{thm:count}]
    Let us first recall the setup.
    Let $v_1,\dots,v_k \in \ZZ^2$ be non-parallel primitive vectors, fix $c_1,\dots,c_k \in \NN$
    and consider a zonotope $Z = \sum_{i=1}^k [-c_i,c_i] v_i$
    that is the Minkowski sum of segments in these directions.
    Let $c>0$ and $n$ be sufficiently large.
    Let $\Omega_n = nZ \cap \ZZ^2$ and let $A_n$ be the incidence matrix
    with rows indexed by the set $L$ of all lines parallel to some $v_i$ that meet $\Omega_n$,
    i.e.~$(A_n)_{\ell,x} = 1_{x \in \ell}$ for each $\ell \in L$ and $x \in \Omega_n$.
    Fix $b = A_n x^*$ for some $x^*\in\{0,1\}^{\Omega_n}$.
    Let $N(b)$ count $x \in \{0,1\}^{\Omega_n}$ with $A_n x = b$.
    Let $p \in [0,1]^{\Omega_n}$ be a vector with $A_n p = b$ maximising the binary entropy
    given by
    \[
        H(p):=-\sum_{v\in\Omega_n}(p_v\log_2p_v+(1-p_v)\log_2(1-p_v)).
    \]
    Suppose $p_v \in [c,1-c]$ for all $v \in \Omega_n$.
    We will show that $N(b) = 2^{H(p)+o(n^{1.51})}$.

    For the upper bound, consider the uniform distribution $U$ on $x \in \{0,1\}^{\Omega_n}$ with $A_n x = b$.
    Let $u\in [0,1]^{\Omega_n}$ be given by $u_v=\PP(U_v=1)$ for each $v\in \Omega_n$.
    Then $A_n u=\EE[A_n U]=b$, so by the maximality of $p$, we have $H(u)\leq H(p)$.
    Thus
    \[
        \log_2 N(b)=H(U)\leq \sum_{v\in\Omega_n}H(U_v)=H(u)\leq H(p),
    \]
    so $N(b) \le 2^{H(p)}$.

    For the lower bound, sample $y\in \{0,1\}^{\Omega_n}$ according to $p$,
    i.e. independently set $y_v=1$ with probability $p_v$ for each $v\in \Omega_n$.
    We will construct a random $x\in\{0,1\}^{\Omega_n}$ with $A_nx=b$ such that $\|x-y\|_1\leq n^{1.51}$ with probability $1-o(n^{-10})$. If the construction fails, we just set $x=x^*$.
    This will suffice. Indeed, since $x-y\in \{0,\pm 1\}^{\Omega_n}$, we have
    \begin{align*}
        H(x-y) & \leq 1+\log_2\paren{\sum_{j=0}^{n^{1.51}} 2^j\binom{|\Omega_n|}{j}}+o(n^{-10})\log_2 3^{|\Omega_n|} \\
               & =O(n^{1.501}\log n)+o(1)=o(n^{1.51}).
    \end{align*}
    Therefore,
    \[\log_2 N(b) \ge H(x) \ge H(y) - H(x-y) \ge H(p)-o(n^{1.51}).\]

    We construct $x$ following the sketch in Subsection~\ref{subsec:absorb}, resolving lines one direction at a time.
    We start with $x=y$ and apply modifications in $k$ stages until $A_n x = b$.
    In stage $i$, we resolve the lines parallel to $v_i$ without disrupting the lines from previous stages.
    That is, $(A_n x)_\ell=b_\ell$ for all lines $\ell$ parallel to $v_j$ with $j\leq i$.

    Each $(A_ny)_\ell$ is a sum of $\Theta(n)$ independent Bernoulli random variables and has expectation $b_\ell$. Chernoff bound and a union bound over $L$ give
    \[
        \PP(\|A_ny-b\|_\infty\geq n^{0.501})
        \leq O(n)\exp(-\Omega(n^{0.002}))
        =\exp(-\Omega(n^{0.002})).
    \]
    Chernoff bounds also show that every line contains $\Omega(n)$ zeroes and $\Omega(n)$ ones in $y$, with probability $1-\exp(-\Omega(n))$.
    To resolve each line $\ell$ parallel to $v_1$, write $n_\ell=b_\ell-(A_nx)_\ell$.
    If $n_\ell>0$, change $n_\ell$ zeroes on $\ell$ to ones.
    If $n_\ell<0$, change $|n_\ell|$ ones on $\ell$ to zeroes.
    We make the $|n_\ell|$ changes one at a time, uniformly at random from the $\Omega(n)$ available options.

    Now suppose $i>1$ and we have resolved all lines in $L$ parallel to any $v_j$ with $j<i$.
    For each $\ell$ parallel to $v_i$ let $n_\ell=b_\ell-(A_nx)_\ell$
    and let $h_i(z)\in\ZZ[z^{\pm}]$ be the generating function with coefficients $n_\ell$.
    We identify $x$ with $f^x\in \ZZ[X^{\pm},Y^{\pm}]$, with $X^s Y^t$ coefficient $x_{s,t}$ for each $(s,t) \in \Omega_n$.
    Similarly we identify $x^*$ with $f^{x^*}\in \ZZ[X^{\pm},Y^{\pm}]$. Writing $v_i=(a_i,b_i)$,
    we have
    \[
        h_i(z)=\pi_{v_i}(f^{x^*}-f^x)\quad\text{where }\pi_{v_i}(f):=f(z^{b_i},z^{-a_i})\in\ZZ[z^{\pm}].
    \]
    As all lines parallel to $v_j$ with $j<i$ are resolved, each $X^{a_j}-Y^{-b_j}$ divides $f^{x^*}-f^x$. Since the $v_j$ are pairwise non-parallel, we have $\prod_{j<i}(X^{a_j}-Y^{-b_j})\mid f^{x^*}-f^x$, and hence $h_i(z)$ is divisible by $\prod_{j<i}(z^{m_{i,j}}-1)$, where $m_{i,j}=|\det(v_i,v_j)|$.
    Since the zonotope $Z$ has diameter $\Theta(n)$, the exponents of $h_i(z)$ are supported in an interval of length $\Theta(n)$.
    Let $m_i=(m_{i,1},\dots,m_{i,i-1})$ and fix $\varepsilon>0$ sufficiently small. By \Cref{lem:decomp-bounded} (after multiplying $h_i$ by a suitable power of $z$ to make it a polynomial with non-zero constant term),
    we can decompose $h_i(z)$ into a sum of at most $C(|h_i|+n)$ (Laurent) polynomials,
    each of which is $m_i$-basic, simple, has diameter at most $\varepsilon n$, and has exponents supported in the interval which supports $h_i(z)$,
    so that every monomial occurs at most $C\|h_i\|_\infty$ times.

    We consider sequentially each $m_i$-basic polynomial in the sum expressing $h_i(z)$.
    When we consider some such polynomial
    \[
        p(z)=\pm z^r\prod_{j<i}(1-z^{m_{i,j}t_j}),\quad r\in\ZZ,\quad t_1,\ldots,t_{i-1}\in\ZZ_+,
    \]
    we modify $x$ by adding some $x^p \in \{-1,0,1\}^{\Omega_n}$ such that $f^{x^p}$ has the form
    \[
        f^{x^p}=\pm X^{a+da_i}Y^{b+db_i}\prod_{j<i}(X^{t_ja_j}-Y^{-t_jb_j}),
    \]
    where $a,b\in\ZZ$ are chosen so that $\pi_{v_i}(f^{x^p})=p(z)$ for every $d\in\ZZ$.
    Since $p$ is simple, all the monomials in $f^{x^p}$ have distinct exponents. In particular, $x^p$ has entries in $\{-1,0,1\}$.
    This modification preserves $n_\ell=0$ for all lines $\ell$ parallel to any $v_j$ with $j<i$.
    Since $p$ has diameter at most $\varepsilon n$ and $\varepsilon$ is sufficiently small,
    there are $\Theta(n)$ choices for $d$ such that $f^{x^p}$ is supported in $\Omega_n$.
    For this modification to be valid, we need the $1$-entries of $x^p$ to align with zeroes of $x$
    and the $-1$-entries of $x^p$ to align with ones of $x$.
    Since $y_v$ are independent Bernoulli random variables with $\PP(y_v=1)\in[c,1-c]$ for all $v\in \Omega_n$,
    Chernoff bounds give with high probability that there are $\Omega(n)$ valid choices for $d$ with respect to the original $y$.
    At each step, choose uniformly among the translations valid with respect to the current $x$, and implement the corresponding modification. We will show later that there are $\Omega(n)$ valid choices for $d$ with high probability.
    After all polynomials in the decomposition of $h_i(z)$ have been considered, we have resolved all lines parallel to $v_i$ without disrupting the previous lines, and proceed to stage $i+1$, or terminate if we are already at stage $k$.

    If the above randomised algorithm is successful, then it terminates with the required
    $x \in \{0,1\}^{\Omega_n}$ with $A_n x = b$ and $\|x-y\|_1 < \tfrac12 n^{1.51}$. Indeed,
    originally we have $x=y$ and $\| A_n x - b \|_\infty < n^{0.501}$, and so $|A_nx-b|_1=O(n^{1.501})$.
    We had $k=O(1)$ applications
    of \Cref{lem:decomp-bounded}, each increasing the total modification support by an $O(1)$ factor, so in the end the total modification support is $|x-y|=O(n^{1.501})$.
    To see that the algorithm is successful with high probability,
    it suffices to show that for any line $\ell \in L$ only $o(n)$ entries of $\ell$ are modified.
    We will show (by induction) that after stage $i$ at most $C^{2i} n^{0.501} $ entries in any line $\ell$ are modified with high probability. Observe that each stage uses only $O(n^{1.501})$ basic polynomials, and so modifies only $O(n^{1.501})$ entries of $x$.

    Consider a line $\ell$ parallel to $v_j$. Consider a modification at stage $i$.
    If $j=i$ then the required bound holds as any monomial occurs at most
    $C\|h_i\|_\infty \le C^{2i-1} n^{0.501}$ times.
    If $j \ne i$, then in stage $i$ we choose at most $C^{2i-1} n^{1.501}$ polynomials one at a time.
    Each polynomial has $\Theta(n)$ choices chosen uniformly at random, at most $2^k$ of which has a monomial on $\ell$.
    Thus the probability that a given polynomial has a monomial on $\ell$ is $O(1/n)$.
    Therefore, the number of polynomials in stage $i$ that have a monomial on $\ell$ is stochastically dominated by a binomial random variable with probability $O(1/n)$ and $C^{2i-1} n^{1.501}$ trials.
    Thus by a Chernoff bound, the probability that more than $C^{2i} n^{0.501}$ polynomials have a monomial on $\ell$ is at most $\exp(-\Theta(n^{0.501}))$.
    A union bound then shows that this holds for all lines $\ell$ simultaneously with high probability.
\end{proof}

\bibliography{main}

\end{document}